\documentclass[oneside,final]{amsart}
\usepackage{geometry}

\usepackage{graphicx}
\usepackage{amsmath}

\usepackage[utf8]{inputenc}
\usepackage{nicefrac}
\usepackage{mathabx}
\usepackage{subcaption}
\usepackage{amsrefs}
\usepackage{esint}
\usepackage{mathtools}
\usepackage{graphicx}
\usepackage{enumerate}

\usepackage{textcomp} 

\usepackage{wrapfig}
\usepackage{xcolor}
\usepackage{enumitem}
\usepackage{verbatim}
\usepackage{hyperref}
\usepackage[right]{showlabels}

\usepackage{cleveref}
\usepackage{subcaption}
\usepackage{autonum}
\usepackage[normalem]{ulem}

\usepackage{enumerate}

\newtheorem{theorem}{Theorem}[section]
\newtheorem{lemma}[theorem]{Lemma}

\theoremstyle{definition}
\newtheorem{definition}[theorem]{Definition}
\newtheorem{example}[theorem]{Example}

\theoremstyle{plain}
\newtheorem{cor}[theorem]{Corollary}

\allowdisplaybreaks
\theoremstyle{remark}
\newtheorem{remark}[theorem]{Remark}

\allowdisplaybreaks
\newcommand{\Span}{\mathrm{span}}

\newcommand{\dist}{\mathrm{dist}}

\newcommand{\mres}{\mathbin{\vrule height 1.6ex depth 0pt width
0.13ex\vrule height 0.13ex depth 0pt width 1.3ex}}

\makeatletter
\newcommand*\bigcdot{\mathpalette\bigcdot@{.5}}
\newcommand*\bigcdot@[2]{\mathbin{\vcenter{\hbox{\scalebox{#2}{$\m@th#1\bullet$}}}}}
\makeatother

\newcommand{\R}{\mathbb{R}}
\newcommand{\RR}{\mathbb{R}}
\newcommand{\N}{\mathbb{N}}
\newcommand{\Z}{\mathbb{Z}}
\newcommand{\W}{\mathcal{W}}
\newcommand{\Div}{\operatorname{div}}
\renewcommand{\H}{\mathcal{H}}

\newcommand{\spt}{\mathrm{spt}}

\makeatletter
\@namedef{subjclassname@2020}{%
  \textup{2020} Mathematics Subject Classification}
\makeatother

\begin{document}
\title{Varifolds with small Willmore energy are weak immersions}
\title{Topology of compact integral 2-varifolds with square integrable mean curvature}
\title{Integral 2-varifolds and $W^{2,2}$ Lipschitz immersions}
\title[Global regularity of integral 2-varifolds with $L^2$ mean curvature]{Global regularity of integral 2-varifolds with square integrable mean curvature}

\author[F.~Rupp]{Fabian Rupp}
\address[F.~Rupp]{Faculty of Mathematics, University of Vienna, Oskar-Morgenstern-Platz 1, 1090 Vienna, Austria.}
\email{fabian.rupp@univie.ac.at}
\author[C.~Scharrer]{Christian Scharrer}
\address[C.~Scharrer]{Institute for Applied Mathematics, University of Bonn, Endenicher Allee 60, 53115 Bonn, Germany.}
\email{scharrer@iam.uni-bonn.de}




\date{\today}

\begin{abstract}
We provide sharp sufficient criteria for an integral $2$-varifold to be
induced by a $W^{2,2}$-conformal immersion of a smooth surface. Our approach is based on a fine analysis of the Hausdorff density  for $2$-varifolds with critical integrability of the mean curvature and a recent local regularity result by Bi--Zhou.
In codimension one, there are only three possible density values below $2$, each of which can be attained with equality in the Li--Yau inequality for the Willmore functional by the unit sphere, the double bubble, and the triple bubble. 
We show that below an optimal 
threshold for the Willmore energy, a varifold induced by a current without boundary is in fact a curvature varifold with a uniform bound on its second fundamental form. 
Consequently, the minimization of the Willmore functional in the class of curvature varifolds with prescribed even Euler characteristic provides smooth solutions for the Willmore problem.
In particular, the ``ambient'' varifold approach and the ``parametric'' approach are equivalent for minimizing the Willmore energy.

\end{abstract}
\maketitle

\section{Introduction}
The Willmore functional is a conformal invariant which quantifies the bending of an immersed surface $F\colon \Sigma\to \R^n$ by 
\begin{align}
    \mathcal{W}(\Sigma) \vcentcolon = \frac{1}{4}\int_\Sigma|H|^2\mathrm{d}\mu,
\end{align}
where $H$ is the 
trace of the second fundamental form
and $\mu$ the surface measure induced by the immersion $F$. While already studied by Blaschke and Thomsen in the 1920's, the functional became more popular through the work of Willmore in the 1960's \cite{Willmore65}. Among closed surfaces, it is not difficult to see that $\mathcal{W}(\Sigma)\geq 4\pi$ with equality only for round spheres \cite[Theorem 7.2.2]{Willmore_Riemannian}. However, 
the minimizer among tori in $\R^3$, the Clifford torus, has only been found rather recently \cite{MarquesNeves14}.
For all $p\in \N_0$ and $n\geq 3$, the geometric variational problem
\begin{align}\label{eq:Willmore_problem}
    \beta_p^n \vcentcolon = \min\{ \mathcal{W}(F) : F\in C^\infty(\Sigma;\R^n) \text{ immersion, $\Sigma$ oriented, }\operatorname{genus}\Sigma=p\},
\end{align}
sometimes referred to as \emph{Willmore problem}, admits a solution \cite{Simon,BauerKuwert} with $\beta_p^n <8\pi$ by \cite{Kusner}. The precise value of $\beta_p^n$ 
is unknown, in general. Besides its large invariance group, a major challenge in the variational analysis of the Willmore energy is that it is Sobolev-critical, thus requiring the development of new techniques to apply the direct method. Historically, two successful approaches have been established to face these issues. As in the pioneering work of Simon \cite{Simon}, immersed surfaces can be viewed as subsets of Euclidean space, see \cite{MR1900754,MR1882663,MR2119722}. This so-called \emph{ambient approach} fits well with interpreting surfaces as measures and allows for applying methods from Geometric Measure Theory, see also \cite{MarquesNeves14}.
In this context, a natural space with good compactness properties is given by the class of integral $2$-varifolds with mean curvature in $L^2$, see \cite{Allard}. Secondly, the \emph{parametric approach} studies surfaces by considering their immersions as maps in a suitable Sobolev space, the space of $W^{2,2}$-conformal Lipschitz immersions \cite{MR2928715,MR2430975,MR3276154}. In accordance with typical bubbling phenomena for Sobolev-critical problems in Geometric Analysis, weak limits may in general develop singularities. Nevertheless, these can be excluded along a minimizing sequence for \eqref{eq:Willmore_problem},see \cite[Theorem 5.3]{MR2928715} and \cite[Theorem 1.1]{Tristan}, by considering suitable comparison surfaces \cite{Kusner,BauerKuwert}, yielding smoothness of minimizers \cite{MR2430975}.

Consequently, the minimization problem \eqref{eq:Willmore_problem} can equivalently be considered in the class $\mathcal{E}_\Sigma$ of $W^{2,2}$-conformal Lipschitz immersions where $\Sigma$ is the unique orientable closed surface of genus $p\in \N_0$. This is a major 
improvement when studying \eqref{eq:Willmore_problem} compared to
the fully ambient approach, i.e., minimization in the class of varifolds, which do not easily allow for a notion of genus. Nonetheless, also among varifolds (having no boundary, in a suitable sense), the global Willmore minimizers are known to be the round spheres as a consequence of the monotonicity formula \cite{MR2119722}.

Integral varifolds comprise the weakest class of generalized surfaces to study the Willmore functional on and, in particular, every $W^{2,2}$-conformal Lipschitz immersion induces an integral curvature varifold, cf.\ \cite[Section 2.2]{MR2928715}. In this article, we prove that, below certain sharp energy thresholds, all integral varifolds can be obtained in this fashion.

\begin{theorem}\label{thm:main_6pi}
Let $\mu$ be an integral 2-varifold in $\R^n$. If $\mu(\R^n)<\infty$ and $\W(\mu)<6\pi$, then $\mu$ is induced by a conformal embedding $F\in \mathcal{E}_\Sigma$ of a smooth compact and connected surface $(\Sigma,g_0)$ without boundary.

The $6\pi$-threshold is sharp: There exists an integral $2$-varifold $\mu_0$ in $\R^3$ with $\mu_0(\R^3)<\infty$ and $\mathcal W(\mu_0)=6\pi$ which satisfies $\Theta^2(\mu_0,x_0)=3/2$ for some $x_0\in\R^3$. In particular, $\mu_0$ is not induced by a Lipschitz immersion.
\end{theorem}

In the case of codimension one, one may employ \cite{MarquesNeves14} to even conclude $\Sigma =\mathbb{S}^2$ in \Cref{thm:main_6pi}. In this particular case, \Cref{thm:main_6pi} is proven independently in a recent preprint by Bi--Zhou \cite{bi2024optimal}.
The energy condition can be weakened 
for integral curvature varifolds \cite{MR0825628,MR1412686}, i.e., integral varifolds with a weak notion of second fundamental form. 

\begin{theorem}\label{thm::main_8pi_curv}
    Let $\mu$ be an integral curvature varifold in $\R^n$ with $\mu(\R^n)<\infty$ and second fundamental form $B\in L^2(\mu)$. If $\W(\mu)<8\pi$, then $\mu$ is induced by a conformal embedding $F\in \mathcal{E}_\Sigma$ of a smooth compact and connected surface $(\Sigma,g_0)$ without boundary. For $n=3$, $\Sigma$ is orientable.
    
     The $8\pi$-threshold is sharp: For all $\varepsilon>0$ there exists an integral curvature $2$-varifold $\mu_0$ in $\R^3$ with $\mu_0(\R^3)<\infty$, $B_{\mu_0}\in L^\infty(\mu_0)$, and $\mathcal W(\mu_0)\le 8\pi+\varepsilon$ which is not induced by a $W^{2,2}$-conformal Lipschitz immersion.
\end{theorem}

The sharpness of the $6\pi$-threshold in \Cref{thm:main_6pi} can be seen by studying a standard double bubble (see \Cref{sec:double-bubble}), an important example from the theory of minimal surfaces \cite{MR1906593}. The fact that this example separates two volumes in $\R^3$ with no reasonable way to orient the separating interface can be seen as a topological obstruction to being induced by an embedding. Such an issue can be avoided by considering varifolds that are the mass measures of $2$-currents without boundary. In fact, below the threshold in \Cref{thm::main_8pi_curv}, these varifolds are already integral curvature varifolds and orientable, in a suitable sense.
\pagebreak
\begin{theorem}\label{cor:KLS}
    For every $\delta>0$ and $n\in\N$, there exists $C(\delta,n)\in (0,\infty)$ with the following property.
      Let $\mu$ be an integral $2$-varifold in $\R^n$ induced by an integral $2$-current $T$ in $\RR^n$ with $\partial T=0$. If $\mu(\RR^n)<\infty$ and $\W(\mu)\leq 8\pi-\delta$, then $\mu$ is an integral curvature varifold with
      \begin{align}
          \Vert B\Vert_{L^2(\mu)} \leq C(\delta,n)
      \end{align}
      and $\spt\,\mu$ is orientable with Euler characteristic $\chi(\spt\,\mu) = 2-2p$ for some $p\in\N_0$.
\end{theorem}

The assumption on the varifold $\mu$ being induced by a current without boundary is naturally compatible with variational settings under (Hausdorff)  density bounds. Indeed, such varifolds are contained in the class of \emph{volume varifolds} that have recently been studied in \cite{ChrisFabian,scharrer2023properties} in order to model cell membranes. Volume varifolds enjoy strong compactness properties and satisfy the assumptions of \Cref{cor:KLS} provided they have unit density almost everywhere \cite{scharrer2023properties}. 

Our strategy to prove \Cref{thm:main_6pi,thm::main_8pi_curv,cor:KLS} is based on a detailed analysis of the density and its relation with the Willmore energy in terms of the Li--Yau inequality \cite{LY,MR2119722}. By \cite{Allard}, the varifold blow up at any fixed point $x_0\in\spt\,\mu$ is a stationary $2$-cone $C$ satisfying $\Theta^2(C,0)=\Theta^2(\mu,x_0)$. Its intersection with the unit sphere is a stationary $1$-varifold $\gamma$ whose total length is given by the cone's density: $\gamma(\mathbb S^{n-1})=2\pi\Theta^2(C,0)$. Any such stationary $1$-varifold is a geodesic net~\cite{AllardAlmgren}, satisfying $\Theta^1(\gamma,p)=\Theta^2(C,p)\le \Theta^2(C,0)$ for $p\in\spt\,\gamma$. In the unit $2$-sphere, isogonal geodesic nets have been fully classified by Heppes \cite{Heppes}. As a byproduct of our density analysis, we find that all these nets with length less than $4\pi$ can be realized with equality in the Li--Yau inequality, see \Cref{lem:density_6pi}.
Under the assumptions of \Cref{thm:main_6pi,thm::main_8pi_curv,cor:KLS}, we can prove that the density of $\mu$ is constantly one in the support of the varifold. We may then apply the recently proven critical case of Allard's regularity theorem by Bi--Zhou \cite{BiZhou22} at \emph{every point,} yielding a local $W^{2,2}$-conformal Lipschitz parametrization of the varifold. Last, we employ the theory of conformal mappings to show that these parametrizations determine a smooth structure.

Below the respective energy thresholds, \Cref{thm:main_6pi,thm::main_8pi_curv,cor:KLS} allow us to extend various concepts and results from the theory of $W^{2,2}$-conformal Lipschitz immersions to varifolds. Concerning regularity theory, using Rivi\`ere's formulation of the Euler--Lagrange equation \cite{MR2430975}, we prove a smoothness result, \Cref{thm:boundary_reg}, for Willmore-minimizing varifolds with prescribed boundary that were introduced in \cite{MR4141858}. 
On the topological side, we define a weak notion of Euler characteristic for integral curvature varifolds  \cite{MR0825628,MR1412686} by means of the classical Gauss--Bonnet theorem and show in \Cref{cor:Gauss_Bonnet} that this is well-defined and an integer under the assumptions of \Cref{thm:main_6pi,thm::main_8pi_curv,cor:KLS}.  
In the context of the Willmore problem \eqref{eq:Willmore_problem}, we can thus show that both established approaches introduced in \cite{Simon} and \cite{MR2928715,MR3276154} are equally as powerful and equivalent from the perspective of the Calculus of Variations.

\begin{cor}[Equivalent approaches for the Willmore problem]\label{cor:equivalent_minimization}
    Let $p\in \N_0, n\geq 3$. Then the minimizers for $\mathcal{W}$ in the following classes coincide.
    \begin{enumerate}[label=(\alph*)]
        \item\label{item:min_W_smooth} The set of smooth immersions $F\colon \Sigma\to\R^n$ where $\Sigma$ is a closed orientable surface with $\operatorname{genus}\Sigma=p$.
        \item\label{item:min_W_conf} The set of $W^{2,2}$-conformal Lipschitz immersions $F\in \mathcal{E}_\Sigma$ into $\R^n$ for $\Sigma$ a closed Riemann surface with $\operatorname{genus}\Sigma=p$.
        \item\label{item:min_W_vari} The set of integral curvature varifolds $\mu$ in $\R^n$ with $\mu(\R^n)<\infty$ and second fundamental form $B\in L^2(\mu)$ such that $\spt\,\mu$ has Euler characteristic $2-2p$.
        \item\label{item:min_W_vari_current} The set of integral $2$-varifolds $\mu$ in $\R^n$ with $\mu(\R^n)+\mathcal{W}(\mu)<\infty$ induced by an integral $2$-current without boundary such that $\spt\,\mu$ has Euler characteristic $2-2p$.
    \end{enumerate}
\end{cor}

Thus, in the class of (curvature) varifolds, we not only have existence  but also the absence of Lavrentiev's phenomenon \cite{MR1553097} for the Willmore problem  \eqref{eq:Willmore_problem}. The latter observation might be relevant for approximating minimizers of $\mathcal{W}$ numerically using nonsmooth shapes. In addition, in \Cref{subsec:approx}, we prove that varifolds satisfying the assumptions of \Cref{thm:main_6pi} or \Cref{thm::main_8pi_curv} can be approximated strongly by smooth surfaces and discuss this in the context of the general approximation problem, see \cite{MR1704565}. Lastly, in \Cref{subsec:Helfrich} we present an extension of \Cref{cor:KLS} to the Helfrich functional which appears in the modelling of lipid bilayers \cite{Helfrich} and also allows for a multiplicity control, see \cite{ChrisFabian}.

We now briefly outline the structure of this article. We first review some necessary background and notation in \Cref{sec:prelims}. In \Cref{sec:density}, we then lay the foundation for proving \Cref{thm:main_6pi,thm::main_8pi_curv,cor:KLS} by discussing density control and its consequences. 
This allows us to complete the proofs of the affirmative statements in the main results, \Cref{thm:main_6pi,thm::main_8pi_curv,cor:KLS}, in \Cref{sec:proof}.
In \Cref{sec:applications}, we then discuss the various applications mentioned in the above paragraphs. \Cref{sec:examples} is devoted to provide several important examples to illustrate the sharpness and significance of our results: the double and triple bubble (\Cref{sec:double-bubble,sec:triple-bubble}), a branched surface (\Cref{sec:branched_immersion}), a surface with an arbitrarily large singular set (\Cref{sec:singular_set}), and an example for an admissible boundary datum for the regularity result, \Cref{thm:boundary_reg}(\Cref{subsec:ex_circ_boundary}). 

\section{Preliminaries}\label{sec:prelims}

Let $U$ be an open subset of $\R^n$, $k\leq n$ be a positive integer, and $\mu$ be a Radon measure over $U$. The \emph{support} of $\mu$ is defined by
\begin{equation}
    \spt\,\mu\vcentcolon=U\setminus\{x\in U\colon \mu(B_r(x))=0\text{ for some }r>0\},
\end{equation}
where $B_r(x)\vcentcolon=\{y\in \R^n\colon |x-y|\le r\}$ are the closed balls in $U$. The \emph{density} at $x\in U$ is defined by
\begin{equation}
    \Theta^k(\mu,x)\vcentcolon=\lim_{r\to0+}\frac{\mu(B_r(x))}{\omega_kr^k}
\end{equation}
provided the limit exists and is finite, where $\omega_k=\mathcal L^k(B_1(0))$ is the Lebesgue measure of the unit ball in $\R^k$. 

An integral $k$-varifold $\mu$ in $U$  is a Radon measure over $U$ of the form
$\mu = \theta \H^k \mres M$, where $\theta\in L^1_{\mathrm{loc}}(M;\N)$ and $M=\{\theta>0\}\subset U$ is an $\mathcal{H}^k$-rectifiable set, that is $M$ is $\mathcal H^k$-measurable and, up to an $\mathcal H^k$-null set, can be covered by countably many $k$-dimensional $C^1$-submanifolds of $\R^n$. A $k$-dimensional linear subspace of $\R^n$ is called \emph{approximate tangent space} of $\mu$ at $x_0\in U$ and denoted with $T_{x_0}\mu$ if 
\begin{equation}
    \lim_{\lambda\to0+}\int_{\lambda^{-1}(M-x_0)}\varphi(x)\theta(x_0+\lambda x)\,\mathrm d\mathcal H^k(x) = \theta(x_0)\int_{T_{x_0}\mu}\varphi(x)\,\mathrm d\mathcal H^k(x)\qquad\text{for $\varphi\in C^0_c(\R^n)$}.
\end{equation}
This exists for $\mu$-almost all $x_0\in U$ \cite[Theorem 11.6]{Simon83}.
The tangential divergence is given by
\begin{equation}
    (\Div_\mu \Phi)(x) \vcentcolon= \sum_{i=1}^k\langle \mathrm D\Phi(x)b_i,b_i\rangle
\end{equation}
whenever $T_x\mu$ exists and $\{b_1,\ldots,b_k\}$ is an orthonormal basis of $T_x\mu$, where $\langle\cdot,\cdot\rangle$ denotes the Euclidean inner product.
We say that $\mu$ has \emph{(generalized) mean curvature} in $U$ if there exists $H\in L^1_{\mathrm{loc}}(\mu;\R^n)$ such that we have the \emph{first variation identity}
\begin{align}\label{eq:first_vari}
    \int \Div_\mu \Phi\,\mathrm{d}\mu  = - \int \langle \Phi, H\rangle\, \mathrm{d}\mu\qquad\text{for $\Phi\in C_c^1(U;\R^n)$}.
\end{align}
For such varifolds, the \emph{Willmore energy} is given by the $L^2$-norm
\begin{align}
    \W(\mu) \coloneqq \frac{1}{4}\int_U |H|^2\mathrm{d}\mu.
\end{align}
Moreover, $\mu$ is called \emph{stationary} in $U$ if and only if $\mu$ has vanishing mean curvature $H\equiv0$.

There also exists a weak notion of second fundamental form for varifolds \cite{MR0825628,MR1412686}. For an integral $2$-varifold $\mu$ in an open subset $U$ of $\R^n$, we denote by $P(x) \in \RR^{n\times n}$ the orthogonal projection onto the approximate tangent space $T_x\mu$ for $\mu$-almost every $x\in U$. We follow the coordinate free notation in \cite{kuwert2023curvature} and term an integral $2$-varifold $\mu$ in $\R^n$ an \emph{integral curvature varifold,} if
there exists a function $B\in L^1_{\mathrm{loc}}(\mu)$ with $B(x)\in BL(\RR^n\times \RR^n;\RR^n)$ for $\mu$-a.e.\ $x\in U$ such that
\begin{align}
    \int \Big(\mathrm{D}_P \Phi(x,P(x)) \cdot B(x) + \langle \mathrm{tr}\, B(x),\Phi(x,P(x))\rangle + \langle \mathrm{D}_x\Phi(x,P(x)), P(x)\rangle \Big)\, \mathrm{d}\mu(x) = 0
\end{align}
for all $\Phi = \Phi(x,P)\in C_c^1(U\times \RR^{n\times n};\RR^n)$.

Let $\bigwedge_k\R^n$ be the space of $k$-vectors and $\bigwedge^k\R^n$ be its dual. An integral $k$-current $T$ in $U$ is a continuous linear functional on $C_c^\infty(U,\bigwedge^k\R^n)$ of the form
\begin{equation}\label{eq:current}
    T(\omega)=\int_{M}\langle \omega(x),\xi(x)\rangle\,\mathrm d\mathcal H^k(x)
\end{equation}
where $M$ is an $\mathcal H^k$-rectifiable subset of $\R^n$, and $\xi$ is an $\mathcal H^k\mres M$-measurable function taking values in $\bigwedge_k\R^n$ such that for $\mathcal H^k\mres M$-almost all $x\in U$, we have $|\xi(x)|\in\N$ and $\xi(x)$ \emph{orients} the approximate tangent space $T_x\mu$ of the induced measure $\mu=|\xi|\mathcal H^k\mres M$, that is $T_x\mu=\{v\in\R^n\colon \xi\wedge v=0\}$. Its boundary is the $(k-1)$-current
\begin{equation}
    \partial T(\omega)\vcentcolon=T(\mathrm d\omega)
\end{equation}
where $\mathrm d$ denotes the exterior derivative. 
In particular, each integral $k$-current $T$ as in $\eqref{eq:current}$ induces an integral $k$-varifold $\mu = \theta \mathcal H^k\mres M$ with $\theta\vcentcolon=|\xi|$.

Given a smooth surface $\Sigma$ with a smooth reference metric $g_0$, a map $F\in W^{2,2}\cap W^{1,\infty}((\Sigma,g_0);\R^n)$ is called a \emph{$W^{2,2}$-Lipschitz immersion} \cite{MR2928715,MR3276154} if there exists $c=c(F,g_0)>0$ with
\begin{align}
    |\mathrm{d}F_p(X)|^2 \geq c |X|^2_{g_0(p)} \quad \text{ for all }X\in T_p\Sigma\label{eq:lipschitz_immersion}
\end{align}
for almost every $p\in \Sigma$. If $\Sigma$ is compact, the definition does not depend on the choice of reference metric $g_0$. We say that $F$ is conformal (with respect to $g_0$) if there exists a function $w\in L^\infty(\Sigma)$, called the \emph{conformal factor} of $F$, such that
\begin{align}
    F^*\langle\cdot,\cdot\rangle  = e^{2w} g_0
\end{align}
almost everywhere in $\Sigma$. The space of $W^{2,2}$-Lipschitz immersions on $\Sigma$ is denoted by $\mathcal{E}_\Sigma$. Any $F\in \mathcal{E}_\Sigma$ induces an integral $2$-varifold with $M=F(\Sigma)$ and $\theta(x) = \# \{p\in \Sigma\mid F(p)=x\}$, see \cite[Section 2.2]{MR2928715}. Similarly, if $\Sigma$ is closed and oriented, any $F\in\mathcal E_\Sigma$ induces an integral $2$-current $T$ in $\R^n$ with $\partial T = 0$. Indeed, we may assume that $\Sigma$ is a submanifold of $\R^m$ and define 
\begin{equation}
    S(\omega)\vcentcolon=\int_{\Sigma}\langle\omega(x),\zeta(x)\rangle\,\mathrm d\mathcal H^2(x),\qquad \text{for $\omega\in C^\infty_c(\R^m,\textstyle\bigwedge^2\R^m)$}
\end{equation}
where in any local positive chart $\varphi$ of $\Sigma$
\begin{equation}
    \zeta = \frac{\partial_1\varphi\wedge \partial_2\varphi}{|\partial_1\varphi\wedge \partial_2\varphi|}. 
\end{equation}
This $2$-current satisfies $\partial S=0$ (see for instance \cite[4.1.31(1)]{Federer}). Now extend $F$ to a Lipschitz map $\bar F\colon \R^m\to\R^n$ and let $T\vcentcolon=\bar F_\# S$ according to \cite[4.1.30]{Federer}. Then $\partial T = \partial \bar F_\#S=\bar F_\#\partial S=0$.

\section{Density analysis}\label{sec:density}

In this section, we will determine the possible values for the Hausdorff density $\Theta^2(\mu,\cdot)$ for varifolds in the energy regime of \Cref{thm:main_6pi,thm::main_8pi_curv,cor:KLS}. We first note the following global regularity result which we derive from \cite{BiZhou22}.

\begin{lemma}\label{lem:induced_surface}
    Let $\mu$ be an integral 2-varifold in $\R^n$ with $\mu(\RR^n)+\W(\mu)<\infty$ and $\Theta^2(\mu,x)=1$ for all $x\in \spt\,\mu$. Then $\mu$ is induced by a conformal embedding $F\in \mathcal{E}_\Sigma$ of a smooth compact surface $(\Sigma,g_0)$ without boundary.
\end{lemma}

\begin{proof}
    The set $\Sigma\coloneq \spt\,\mu$ is compact by the assumptions, see \cite[(A.22)]{MR2119722}.
    Given any $x_0\in \Sigma$, we may take $r=r(x_0)>0$ small enough such that the assumptions of \Cref{thm:BiZhou_Scaled} are satisfied. In particular, by \Cref{thm:BiZhou_Scaled}, we find an open set $U\subset \R^2$ and a conformal bi-Lipschitz homeomorphism $f\colon U\subset \R^2\to f(U)\subset \Sigma$ with $x_0\in f(U)$. Its inverse $f^{-1}$ provides a Lipschitz continuous local chart of $\Sigma$ near $x_0$. By compactness, we may find a finite cover $\Sigma=\bigcup_{\alpha\in I_0} f_\alpha(U_\alpha)$ turning $\Sigma$ into a compact topological $2$-dimensional manifold without boundary, embedded in $\R^n$ by the inclusion. 

    Given two parametrizations $f_\alpha\colon U_\alpha\to f_\alpha(U_\alpha)\subset \Sigma$ and $f_\beta\colon U_\beta\to f_\beta(U_\beta)\subset\Sigma$, we consider the transition map $\sigma \vcentcolon= f_\beta^{-1}\circ f_\alpha \colon  W\vcentcolon= f_\alpha^{-1}(f_\beta(U_\beta)) \to \RR^2$ which is again a bi-Lipschitz homeomorphism onto its image. By Rademacher's theorem, there exist sets $\hat{U}_\alpha\subset U_\alpha$, $\hat{U}_\beta\subset U_\beta$, $\hat{W}\subset W$, each of full measure, such that $f_\alpha,f_\beta,\sigma$ are differentiable on $\hat{U}_\alpha$,  $\hat{U}_\beta$, $\hat{W}$, respectively. If $z\in \hat{U}_\alpha\cap\hat{W}\cap \sigma^{-1}(\hat{U}_\beta)$ all derivatives exist and
    \begin{align}
        \mathrm{D} f_\alpha(z) = \mathrm{D}(f_\beta\circ \sigma)(z) = \mathrm{D}f_\beta(\sigma(z)) \mathrm{D}\sigma(z).
    \end{align}
    Since $\sigma^{-1}$ is Lipschitz, $\hat{U}_\alpha\cap\hat{W}\cap \sigma^{-1}(\hat{U}_\beta)$ is also of full measure in $W$.
    Consequently, for almost every $z\in W$ and $i,j\in \{1,2\}$ we have
    \begin{align}\label{eq:Dxi}
        e^{2w_\beta(\sigma(z))}\langle \partial_i \sigma(z), \partial_j \sigma(z)\rangle = e^{2w_\alpha(z)} \delta_{ij},
    \end{align}
    where $w_\alpha,w_\beta$ are the bounded conformal factors of $f_\alpha,f_\beta$ given by \Cref{thm:BiZhou_Scaled}. Thus $\mathrm{D}\sigma(z)$ is conformal for almost every $z\in W$.
    Since $\sigma$ is injective and Lipschitz, by a degree argument (see \cite[4.1.26]{Federer}) 
    we have that either $\det \mathrm{D}\sigma \geq 0$ a.e.\ in $W$ or $\det \mathrm{D}\sigma \leq 0$ a.e.\ in $W$. In the first case, \eqref{eq:Dxi} implies that $\sigma$ is holomorphic at almost every $z \in W$, whereas in the second case it yields that $\sigma$ is antiholomorphic at almost every $z\in W$. Since $\sigma$ is Lipschitz, in both cases this implies that $\sigma$ is weakly harmonic in $W$ and thus smooth. This gives $\Sigma$ a smooth structure with respect to which the $f_\alpha$ are smooth diffeomorphisms. In particular, the tangent space of $\Sigma$ at $x\in\Sigma$ is given by $\operatorname{span}\{ \partial_1f_\alpha(z),\partial_2f_\alpha(z)\}$ if $x=f_\alpha(z)$.
    
    Choose a smooth partition of unity $\{\eta_\alpha\colon\alpha\in I_0\}$ subordinate to the open cover $\{f_\alpha(U_\alpha)\colon\alpha\in I_0\}$ and let
    \begin{equation}
        g_0\vcentcolon=\sum_{\alpha\in I_0}\eta_\alpha g^\alpha 
    \end{equation}
    where for each $\alpha\in I_0$, $g^\alpha$ is defined on the tangent bundle $Tf_\alpha(U_\alpha)$ by
    $g^\alpha(\partial_if_\alpha,\partial_jf_\alpha) \vcentcolon= \delta_{ij}$.
    This defines a smooth Riemannian metric on $\Sigma$ and, for any $\alpha\in I_0$, we may use \eqref{eq:Dxi} to compute
    \begin{equation}
        g_0(\partial_i f_\alpha,\partial_jf_\alpha) = \sum_{\beta\in I_0}\eta_\beta \langle \partial_i(f_\beta^{-1}\circ f_\alpha),\partial_j(f_\beta^{-1}\circ f_\alpha)\rangle = \sum_{\beta\in I_0}\eta_\beta e^{2(w_\alpha - w_\beta \circ f_\beta^{-1}\circ f_\alpha)} \delta_{ij}.
    \end{equation}
    This implies that the inclusion map $F\colon (\Sigma,g_0)\to \RR^n$ is conformal. Moreover, the $L^\infty$-boundedness of the conformal factors in \Cref{thm:BiZhou_Scaled} and \eqref{eq:Dxi} imply $F\in W^{1,\infty}((\Sigma,g_0);\RR^n)$ and \eqref{eq:lipschitz_immersion}. In a local parametrization $f_\alpha$, we have
    \begin{align}
        \big((\nabla^{g_0})^2_{ij}F\big)\circ f_\alpha = \partial_i \partial_j f_\alpha - \Gamma(g_0)_{ij}^k\partial_k f_\alpha,
    \end{align}
    so that $f_\alpha\in W^{2,2}(U_\alpha;\RR^n)$ yields $F\in W^{2,2}((\Sigma,g_0);\RR^n)$ and thus $F\in \mathcal{E}_{\Sigma}$. Lastly, by construction, we have that $\mu = \mathcal{H}^2\mres \Sigma$ is induced by $F$.
\end{proof}

 Next, we construct a stationary $1$-varifold in $\mathbb S^{n-1}$ from a stationary $2$-cone in $\R^n$. This will allow us to use the classification results from \cite{AllardAlmgren,Heppes} in the sequel.
 
\begin{lemma}\label{lem:cone}
    Let $C$ be a stationary integral $2$-varifold in $\R^n$ such that for all Borel sets $A\subset\R^n$
    \begin{equation}\label{eq:lem:cone}
        C(A) = r^2 C(r^{-1}A)\qquad\text{for all $r>0$}.
    \end{equation}
    Then
    \begin{equation}\label{eq:lem:cone:gamma}
        \gamma(\varphi)\vcentcolon=\int_{\mathbb S^{n-1}}\varphi(p)\Theta^2(C,p)\,\mathrm d\mathcal H^1(p)\qquad \text{for $\varphi\in C^0_c(\R^n)$}
    \end{equation}
    defines a stationary $1$-varifold in $\mathbb S^{n-1}$ with 
    $\Theta^1(\gamma,p)=\Theta^2(C,p)\ge1$ for all $p\in\spt\,\gamma$.
\end{lemma}

\begin{proof}
    In view of the monotonicity formula \Cref{lem:ap:monotonicity}, the support of $C$ is the closed set  
    \begin{equation}
        \Gamma\vcentcolon=\spt\,C=\{x\in\R^n\colon1\le\Theta^2(C,x)<\infty\}
    \end{equation}
    and is $\mathcal H^2$-rectifiable (cf.\ Theorems 2.8(5) and 3.5(1) in \cite{Allard}). Since $B_r(\lambda x)=\lambda B_{\lambda^{-1}r}(x)$ for all $x\in\R^n$ and $r,\lambda>0$, we infer from \eqref{eq:lem:cone} that 
    \begin{equation}\label{eq:lem:cone:density}
        \Theta^2(C,\lambda x) = \Theta^2(C,x).
    \end{equation}
    Moreover, $C(B_r(0))=r^2C(B_1(0))=\pi r^2\Theta^2(C,0)$. Thus, \cite[Theorem 5.2(2)(a)]{Allard} implies
    \begin{equation}\label{eq:lem:cone:tangent_cone}
        T_xC\supset \Span\{x\}\qquad \text{for $\mathcal H^2$-almost all $x\in\Gamma$.}
    \end{equation}
    Let $f(x)\vcentcolon=|x|$ for $x\in\R^n$. Since $f$ is differentiable away from zero, we have that in the notation of \cite[3.2.16]{Federer},
    \begin{equation}\label{eq:lem:cone:approximate_differential}
        (\mathcal H^2\mres\Gamma,2)\,\mathrm{ap}\, \mathrm Df(x) =  \mathrm Df(x)|_{T_xC}\qquad\text{for $\mathcal{H}^2$-almost all $x\in\Gamma$}.
    \end{equation}
    Combining \eqref{eq:lem:cone:tangent_cone} and \eqref{eq:lem:cone:approximate_differential}, we see that in the notation of the coarea formula \cite[3.2.22]{Federer}, it holds
    \begin{equation}
        \mathrm{ap}\, J_1f(x) = \|\textstyle\bigwedge_1(\mathcal H^2\mres\Gamma,2)\,\mathrm{ap}\,\mathrm Df(x)\|=\|\mathrm Df(x)|_{T_xC}\|=1
    \end{equation}
    for $\mathcal H^2$-almost all $x\in \Gamma$. It follows 
    \begin{equation}\label{eq:lem:cone:coarea}
        \int_{\Gamma} \varphi \,\mathrm dC= \int_{\Gamma}\varphi(x)\Theta^2(C,x)\,\mathrm{ap}\, J_1f(x)\,\mathrm d\mathcal H^2(x) =\int_0^\infty\int_{r\mathbb S^{n-1}}\varphi(p)\Theta^2(C,p)\,\mathrm d\mathcal H^1(p)\,\mathrm dr
    \end{equation}
    for $\varphi\in C^0_c(\R^n)$. Thus, we can combine Lemma 2.6(3) and Theorem 5.2(2) of \cite{Allard} with \eqref{eq:lem:cone:coarea} to deduce that for $\mathcal L^1$-almost all $r>0$, 
    \begin{equation}
        B(r)(\varphi)\vcentcolon=r^{-1}\int_{r\mathbb S^{n-1}}\varphi(r^{-1}p)\Theta^2(C,p)\,\mathrm d \mathcal H^1(p)
    \end{equation}
    is a stationary $1$-varifold in $\mathbb S^{n-1}$. By \eqref{eq:lem:cone:density}, $B(r)=\gamma$ is $\mathcal L^1$-almost constant. Thus, the conclusion follows from \cite[Theorem 5.2(2)(f)]{Allard}.
\end{proof}

\begin{lemma}\label{lem:density_6pi}
    Let $U\subset\R^n$ be open. Suppose $\mu$ is an integral $2$-varifold in $U$ with $\mathcal W(\mu) <\infty$. 
    Then the following hold.
    \begin{enumerate}
        \item \label{it:lem:6pi} If $a\in U$ with $1\le\Theta^2(\mu,a)<3/2$, then $\Theta^2(\mu,a)=1$.
        \item \label{it:lem:6pi2} If $n=3$ and $a\in U$ with $1\le\Theta^2(\mu,a)<2$, then $\Theta^2(\mu,a)\in\{1,3/2,3\arccos(-1/3)/\pi\}$. Moreover, for $U=\R^3$, all values can be obtained with equality in the Li--Yau inequality:
        \begin{enumerate}
            \item \label{it:6pi:sphere} $\Theta^2(\mu,x_0)=1=\frac{1}{4\pi}\mathcal W(\mu)$ for $\mu$ the unit sphere and any $x_0\in\spt\,\mu$;
            \item \label{it:6pi:double} $\Theta^2(\mu,x_0)=\frac{3}{2}=\frac{1}{4\pi}\mathcal W(\mu)$ for $\mu$ the standard double bubble and $x_0$ on a circle;
            \item \label{it:6pi:triple} $\Theta^2(\mu,x_{1/2})=\frac{3\arccos(-\frac{1}3)}{\pi}=\frac{1}{4\pi}\mathcal W(\mu)$ for $\mu$ a triple bubble and precisely two points $x_{1/2}$.
        \end{enumerate}
    \end{enumerate}
\end{lemma}

\begin{proof}[Proof of \Cref{lem:density_6pi}]
    Let $a\in\spt\,\mu$ be given. 
    We apply \cite[3.4(1)(2)]{Allard} to obtain a sequence $r_k$ of positive numbers diverging to infinity such that
    \begin{equation}
        C(\varphi)\vcentcolon=\lim_{k\to\infty}r_k^2\int \varphi\bigl(r_k(x-a)\bigr)\,\mathrm d\mu(x)\qquad\text{for $\varphi\in C_c^0(\R^n)$}
    \end{equation}
    defines a Radon measure over $\R^n$ with $\Theta^2(C,0)=\Theta^2(\mu,a)$. Moreover, by \cite[6.5, 5.2(2)(b)]{Allard}, $C$ satisfies the assumptions of \Cref{lem:cone}. Therefore, the structure theorem for stationary $1$-varifolds \cite{AllardAlmgren} can be applied to 
    \begin{equation}
         \gamma(\varphi)\vcentcolon=\int_{\mathbb S^{n-1}}\varphi(p)\Theta^2(C,p)\,\mathrm d\mathcal H^1(p)\qquad \text{for $\varphi\in C^0_c(\R^n)$}
    \end{equation}
    as a stationary $1$-varifold in $\mathbb S^{n-1}$. It follows that $\gamma$ consists of geodesic arcs and, since $\Theta^2(C,p)<\infty$ for all $p\in \spt\,\gamma$ by \Cref{lem:ap:monotonicity}, there exists an integer $k\ge2$ such that 
    \begin{equation}\label{eq:lem:density:gamma}
        \Theta^2(C,p)=\Theta^1(\gamma,p)=k/2.
    \end{equation}
    If this equation is true for some $p\in \spt\,\gamma$ with $k>2$, then $p$ is a junction of geodesic arcs whose multiplicity adds up to $k$ and Equation \eqref{eq:lem:cone:density} combined with the upper semi-continuity of $\Theta^2(C,\cdot)$ implies $\Theta^2(C,0)\ge 3/2$. In order to prove the first statement, we may thus assume $k=2$ for all $p\in \spt\,\gamma$. Since $\gamma$ is stationary, the $2$ geodesic arcs meeting at $p$ can only meet tangentially. Therefore, $\gamma$ is a union of $m\in \N$ great circles. Hence, $C$ is a union of $m$ multiplicity-$1$ planes meeting in the origin. Since $\Theta^2(C,0)=\Theta^2(\mu,a) < 3/2$, we have $m=1$, implying that $C$ is a multiplicity-$1$ plane, and $\Theta^2(\mu,a)=1$. 
    
    Combining \eqref{eq:lem:cone:density} and \eqref{eq:lem:cone:coarea}, we see that 
    \begin{align}\label{eq:theta_vs_length}
        \Theta^2(C,0)=\frac{\gamma(\mathbb S^{n-1})}{2\pi}.    
    \end{align}
    Hence, the first part of Statement \eqref{it:lem:6pi2} follows from \Cref{lem:ap:Heppes}. Statement \eqref{it:6pi:sphere} is obvious whereas \eqref{it:6pi:double} and \eqref{it:6pi:triple} follow from Sections \ref{sec:double-bubble} and \ref{sec:triple-bubble}.
\end{proof}

Next, we show that unit density is the only possibility under the condition to be induced by a current as in \Cref{cor:KLS}. 

\begin{lemma} \label{lem:main_8pi}
    Let $U\subset\R^n$ be open. Suppose $\mu$ is the varifold induced by an integral $2$-current in $U$ with $\partial T = 0$. If $\Theta^2(\mu,x)<2$ for all $x\in\spt\,\mu$ and $\mathcal W(\mu) <\infty$, then there holds $\Theta^2(\mu,x)=1$ for all $x\in\spt\,\mu$.
\end{lemma}

\begin{proof}
    The current $T$ can be represented as
    \begin{equation}\label{eq:thm:8pi:current}
        T(\omega) = \int_{\Sigma}\langle\omega,\xi\rangle\,\mathrm d \mathcal H^2 \qquad\text{for $\omega\in C^\infty_c(U,\textstyle\bigwedge^2\R^n)$}
    \end{equation}
    where $\Sigma=\spt\,\mu$ and $\xi(x)$ is a simple $2$-vector orienting $T_x\mu$ with $|\xi(x)|=1$ for $\mathcal H^2$-almost all $x\in \Sigma$. Let $x_0\in\Sigma$. 
    By \cite[3.4(1)(2), 4.12(2), 6.5]{Allard}, there exists a sequence $r_k$ such that $\lim_{k\to\infty}r_k=\infty$ and
    \begin{equation}
        C(\varphi)\vcentcolon=\lim_{k\to\infty}r_k^2\int_{\Sigma} \varphi\bigl(r_k(x-x_0)\bigr)\,\mathrm d\mathcal H^2(x)\qquad\text{for $\varphi\in C_c^0(\R^n)$}
    \end{equation}
    defines a stationary integral $2$-varifold satisfying $\Theta^2(\mu,x_0)=\Theta^2(C,0)$ and satisfying the hypothesis of \Cref{lem:cone}.
    By the monotonicity inequality \Cref{lem:ap:monotonicity},
    \begin{equation}\label{eq:thm:8pi:density}
        1\le \Theta^2(C,x) = \lim_{r\to0}\frac{C(B_r(x))}{\pi r^2}\le \limsup_{s\to\infty}\frac{C(B_s(x))}{\pi s^2} = \Theta^2(C,0)=\Theta^2(\mu,x_0)<2
    \end{equation}
    for all $x\in \spt \, C =\vcentcolon\Gamma$. Define
    \begin{equation}
        \xi_k(x)\vcentcolon=
        \begin{cases}
            \xi(\frac{x}{r_k}+x_0)&\text{for $x\in r_k(\Sigma-x_0)$}\\
            0&\text{for $x\notin r_k(\Sigma-x_0)$}.
        \end{cases}
    \end{equation}
    We apply the theory of slicing by means of the coarea formula, see \cite[Lemma 28.5]{Simon83} and \cite[Theorem 4.10(2)]{Allard} to obtain radii $\rho_k$ diverging to infinity such that the currents
    \begin{equation}
        T_k(\omega)\vcentcolon=\int_{B_{\rho_k}(0)}\langle \omega(x),\xi_k(x)\rangle\,\mathrm d\mathcal H^2(x) \qquad\text{for $\omega\in C^\infty_c(\R^n,\textstyle\bigwedge^2\R^n)$}
    \end{equation}
    and the oriented varifolds
    \begin{equation}\label{eq:thm:8pi:oriented_blow-up}
        V^{\mathrm o}_k(\varphi) \vcentcolon=\int_{B_{\rho_k}(0)}\varphi(x,\xi_k(x))\,\mathrm d\mathcal H^2(x)\qquad \text{for $\varphi \in C^0_c(\R^n\times\mathbb G^\mathrm{o}(n,2))$}.
    \end{equation}
    satisfy
    \begin{equation}
        \sup_{k\in\N}\|\partial T_k\|(K) + \|\delta V_k^\mathrm{o}\|(K) < \infty\qquad \text{for all compact $K\subset \R^n$}.
    \end{equation}
    Here, $\mathbb G^\mathrm{o}(n,2)$ denotes the set of oriented 2-dimensional subspaces in $\R^n$ as defined in \cite{MR0825628}, see also \Cref{subsec:Helfrich}.
    The compactness theorem for oriented integral varifolds \cite[Theorem 3.1]{MR0825628} implies the existence of Borel functions $\nu\colon\Gamma\to\mathbb G^\mathrm{o}(n,2)$ and $\theta_1,\theta_2\colon\Gamma\to\N$ such that
    \begin{equation}\label{eq:thm:8pi:oriented_blow-up_limit}
        \lim_{k\to\infty}V^\mathrm{o}_k(\varphi)=\int_{\Gamma}\Bigl[\varphi(x,\nu(x))\theta_1(x) + \varphi(x,-\nu(x))\theta_2(x)\Bigl]\,\mathrm d\mathcal H^2(x).
    \end{equation}
    By \eqref{eq:thm:8pi:density}, we have
    \begin{equation}\label{eq:thm:8pi:oriented_denisty}
        \theta_1(x)+\theta_2(x)=\Theta^2(C,x)=1\qquad \text{for $\mathcal H^2$-almost all $x\in\Gamma$}.
    \end{equation}
    On the other hand, by the compactness theorem for integral currents (see for instance Theorems 26.14 and 32.2 in \cite{Simon83}), there exists an $\mathcal H^2$-measurable function $\hat\xi$ on $\R^n$ such that $\hat\xi(x)$ is a simple $2$-vector with $|\hat \xi(x)|\in\N$ for $\mathcal H^2$-almost all $x\in\hat\Sigma = \spt\,  \hat T$ and
    \begin{equation}\label{eq:thm:8pi:current_convergence}
        \hat T(\omega)\vcentcolon=\lim_{k\to\infty} T_k(\omega)=\int_{\hat \Sigma} \langle\omega(x),\hat\xi(x)\rangle\,\mathrm d\mathcal H^2(x)\qquad \text{for $\omega\in C^\infty_c(\R^n,\textstyle\bigwedge^2\R^n)$}
    \end{equation}
    defines an integral $2$-current in $\R^n$. From \eqref{eq:thm:8pi:oriented_blow-up} and \eqref{eq:thm:8pi:oriented_blow-up_limit} it follows that
    \begin{equation}
        \hat\xi(x) =\nu(x)\bigl(\theta_1(x)-\theta_2(x)\bigr).
    \end{equation}
    Moreover, \eqref{eq:thm:8pi:oriented_denisty} implies $|\hat\xi|=|\nu|$. Thus, the varifold $C$ is induced by $\hat T$. The weak convergence \eqref{eq:thm:8pi:current_convergence} also implies that $\partial \hat T=0$. 
    
    Let $f\colon \R^n\to\R$ be defined by $f(x)\vcentcolon=|x|$ and $\Gamma_t\vcentcolon=f^{-1}(t)\cap\Gamma$. By \cite[\S 28]{Simon83} there exists $t>0$ and an $\mathcal H^1\mres\Gamma_t$-measurable map $\hat\xi_t\colon\Gamma_t\to\mathbb S^{n-1}$ such that  
    \begin{equation}
        \hat T_t(\omega)\vcentcolon=\int_{\Gamma_t}\langle\omega(p),\hat\xi_t(p)\rangle\Theta^2(C,p)\,\mathrm d\mathcal H^1(p)\qquad \text{for $\omega \in C^\infty_c(\R^n,\textstyle\bigwedge^1\R^n)$}
    \end{equation}
    defines an integral $1$-current with $\partial \hat T_t=0$. From \Cref{lem:cone}  we know that $\hat T_t$ induces a stationary $1$-varifold in $t\mathbb S^{n-1}$. Similarly as in \eqref{eq:lem:orientability:constancy}, we may combine the constancy theorem with the structure theorem \cite{AllardAlmgren} to obtain a finite family of geodesics $\alpha_i\colon [0,1]\to t\mathbb S^{n-1}$ for $i=1,\ldots,N$ such that $\alpha_i|_{[0,1)}$ is injective and
    \begin{equation}
        \hat T_t = \sum_{i=1}^N\sigma_i{\alpha_i}_\#\mathbb E^1\mres [0,1]
    \end{equation}
    where $\sigma_i$ are integers with $|\sigma_i|=1$. Define
    \begin{equation}
        N_i(s)\vcentcolon=\#\Bigl(\{j\colon \alpha_j(0)=\alpha_i(s)\}\cup\{j\colon \alpha_j(1)=\alpha_i(s)\}\Bigr)
    \end{equation}
    From \eqref{eq:thm:8pi:density} we have $2\le N_i(s) <4$ for $s=0,1$. Now, evaluating
    \begin{equation}
        0=\partial\hat T_t = \sum_{i=1}^N\sigma_i(\alpha_i(1)-\alpha_i(0))
    \end{equation}
    at $\varphi\in C_c^\infty(U;\Lambda^0(\R^n)) = C_c^\infty(U)$ with $\varphi$ a cutoff function satisfying $\varphi(\alpha_i(s))=1$,    
    we also know $N_i(s) \neq 3$ for $s=0,1$. Hence, $N_i(s)=2$ for $s=0,1$ and $\hat T_t$ consists of disjoint great circles. By \eqref{eq:thm:8pi:density}, $\hat T_t$ consists of only one great circle, $C$ is a multiplicity-$1$ plane and $\Theta^2(\mu,x_0)=\Theta^2(C,0)=1$.
\end{proof}

Using the classification of curvature varifolds with vanishing second fundamental form \cite{MR0840281}, we can exclude noninteger densities.

\begin{lemma}\label{lem:curv_var}
    Let $U\subset \R^n$ be open and $\mu$ be an integral  curvature $2$-varifold in $U$ with $B\in L^2_\mathrm{loc}(\mu)$. Then $\Theta^2(\mu,a)\in\mathbb N$ for all $a\in \spt\,\mu$.    
\end{lemma}

\begin{proof}
    Let $a\in\spt\,\mu$. By \cite[3.4.(1)(2), 6.5]{Allard} there exists an integral $2$-varifold $C$ in $\R^n$ and a sequence $r_k$ of positive numbers diverging to infinity such that $\Theta^2(\mu,a) =\Theta^2(C,0)$ and such that the $2$-varifolds $\mu_k$ in $r_k(U-a)$ defined by
    \begin{equation}
        \mu_k(\varphi)\vcentcolon=\int_{r_k(U-a)}\varphi(x,P_k(x))\,\mathrm d\mu_k(x) \vcentcolon= r^2_k\int_{U}\varphi(r_k(x-a),P(x))\,\mathrm d\mu(x)
    \end{equation}
    for $\varphi\in C_c(r_k(U-a)\times \R^{n\times n})$ satisfy $P_k(x)\vcentcolon=T_x\mu_k=T_{\frac{x}{r_k}+a}\mu$ for $\mu_k$-almost all $x\in r_k(U-a)$ and
    \begin{equation}
        \lim_{k\to\infty}\mu_k(\varphi)=\int_{\R^n}\varphi(x,T_xC)\,\mathrm dC(x)\qquad \text{for $\varphi\in C_c(\R^n\times \R^{n\times n})$}.
    \end{equation}
    Let $\Phi = \Phi(x,P)\in C_c^1(\R^n\times \RR^{n\times n};\RR^n)$ and define $\Phi_k(x,P)\vcentcolon=\Phi(r_k(x-a),P)$. Then, for $\spt\,\Phi$ small enough or $k$ large enough,
    \begin{align}
        &\int\langle\mathrm D_x\Phi(x,P_k(x)),P_k(x)\rangle\,\mathrm d\mu_k(x) = r_k^2\int \frac{1}{r_k}\langle\mathrm D_x\Phi_k(x,P(x)),P(x)\rangle\,\mathrm d\mu(x) \\
        &\qquad= -r_k^2\int\frac{1}{r_k}\Big(\mathrm{D}_P \Phi_k(x,P(x)) \cdot B(x) + \langle \mathrm{tr}\, B(x),\Phi_k(x,P(x))\rangle\Bigr)\,\mathrm d\mu(x) \\
        &\qquad = -\int\Big(\mathrm{D}_P \Phi(x,P_k(x)) \cdot B_k(x) + \langle \mathrm{tr}\, B_k(x),\Phi(x,P_k(x))\rangle\Bigr)\,\mathrm d\mu_k(x)
    \end{align}
    for $B_k(x)\vcentcolon = \frac{1}{r_k}B(\frac{x}{r_k}+a)$. In particular, $\mu_k$ is an integral curvature varifold and
    \begin{equation}\label{eq:lem:curv_var}
        \int_{B_\rho(0)}|B_k|\,\mathrm d\mu_k \le r_k\int_{B_{\rho/r_k(a)}}|B|\,\mathrm d\mu \le \rho\Bigl(\int_{B_{\rho/r_k(a)}}|B|^2\,\mathrm d\mu\Bigr)^{\frac{1}{2}}\Bigl(\frac{\mu(B_{\rho/r_k}(a))}{(\rho/r_k)^2}\Bigr)^{\frac{1}{2}}.
    \end{equation}
    By \cite[Remark 5.2.3]{MR0825628} (see also \cite[(19)]{kuwert2023curvature}) we have that $\mu$ satisfies \eqref{eq:first_vari} for $H=\mathrm{tr}\, B\in L^2_{\mathrm{loc}}(\mu,\R^n)$. Thus, \Cref{lem:ap:monotonicity} implies $\Theta^2(\mu,a)<\infty$ and the right hand side of \eqref{eq:lem:curv_var} goes to zero as $k\to\infty$. Therefore, by \cite[5.3.2]{MR0825628}, $C$ is an integral curvature varifold and, by weak lower semicontinuity of the $L^1$-norm, we have $B_C=0$. Now, the conclusion follows from \cite[p.\ 292]{MR0840281}.
\end{proof}

As a last key ingredient, we show that being induced by a current as in \Cref{cor:KLS} implies orientability of the surface in \Cref{lem:induced_surface}.

\begin{lemma}\label{lem:orientability}
    Let $\mu$ be the varifold induced by an integral $2$-current $T$ in $\R^n$ with $\partial T=0$. If $\mu(\R^n) + \mathcal W(\mu)<\infty$ and $\Theta^2(\mu,x)=1$ for all $x\in\spt\,\mu$, then the compact surface $(\Sigma,g_0)$ according to \Cref{lem:induced_surface} is orientable, thus a Riemann surface, and $T$ is induced by the embedding $F$ of $\Sigma$. 
\end{lemma}

\begin{proof}
    The current $T$ can be represented as
    \begin{equation}\label{eq:lem:orientability:current}
        T(\omega) = \int_{\spt\,\mu}\langle\omega,\xi\rangle\,\mathrm d \mathcal H^2
    \end{equation}
    where $\xi(x)$ is a simple $2$-vector orienting $T_x\mu$, and $|\xi(x)|=1$ for $\mathcal H^2$-almost all $x\in \spt\,\mu$. Let $x_0\in \Sigma\vcentcolon=\spt\,\mu$. Let $f\colon U\to\R^n$ be a Lipschitz map according to \Cref{thm:BiZhou_Scaled} and $r>0$ such that $B_r(x_0)\cap\Sigma \subset f(U)$. Let $g\colon f(U)\to\R^2$ be the inverse of $f$. Choose Lipschitz continuous extensions $\bar f\colon\R^2\to\R^n$ and $\bar g\colon\R^n\to\R^2$ of $f$ and $g$, respectively. Let $\eta \in C^\infty_c(B_r(x_0))$ such that $\eta(x)=1$ for all $x\in B_{r/2}(x_0)$. Since $\partial T = 0$, we have that 
    \begin{equation}
        \partial (T\mres\eta) = -T\mres \mathrm d\eta 
    \end{equation}
    and hence
    \begin{equation}
        \spt \,\partial \bar g_\#(T\mres\eta)\subset D\setminus g(B_{r/2}(x_0)).
    \end{equation}
    Therefore, in the notation of \cite[4.1.7]{Federer}, the constancy theorem implies the existence of $c\in\R$ such that
    \begin{equation}\label{eq:lem:orientability:constancy}
        \spt\,[\bar g_\#(T\mres\eta) - c(\mathbb E^2 \mres D)] \subset D\setminus g(B_{r/2}(x_0)).
    \end{equation}
    Mapping both currents by $F$, it follows from \cite[Lemma 4.1.25]{Federer} and \eqref{eq:lem:orientability:current} that     \begin{equation}\label{eq:lem:orientability:normal}
        \xi(x) = c\,\frac{\partial_1f\wedge\partial_2f}{|\partial_1f\wedge\partial_2f|}\circ g(x)\qquad\text{for $\mathcal H^2$-almost all $x\in B_{r/2}(x_0)\cap \Sigma$}.
    \end{equation}
    Hence, $|c|=1$ and by swapping the coordinates of $f$, one may achieve $c=1$. 
    Given any two parametrizations $f_\alpha\colon U_\alpha\to\Sigma$, $f_\beta\colon U_\beta\to\Sigma$ as in \Cref{thm:BiZhou_Scaled} satisfying \eqref{eq:lem:orientability:normal} with $c=1$, and denoting the transition map with $\sigma\vcentcolon=f_\beta^{-1}\circ f_\alpha|_{f_\alpha^{-1}(f_\beta(U_\beta))}$, it follows
    \begin{equation}
        |\partial_1f_\alpha \wedge\partial_2f_\alpha|\,\xi = \det(\mathrm D\sigma)(|\partial_1f_\beta\wedge \partial_2f_\beta|\circ\sigma)\,\xi.
    \end{equation}
    Thus, $\det(\mathrm D\sigma)\ge0$, and $\sigma$ is holomorphic, cf.\ the proof of \Cref{lem:induced_surface}. Therefore, the family of parametrizations $f\colon U\to\Sigma$ satisfying \eqref{eq:lem:orientability:normal} with $c=1$ induces a conformal structure on $\Sigma$. In particular, given local coordinates $(z^1,z^2)$ corresponding to a parametrization $f$ of the conformal structure, we see that
    \begin{equation}
        \frac{\partial_{z^1}F\wedge \partial_{z^2}F}{|\partial_{z^1}F\wedge \partial_{z^2}F|} = \frac{\partial_1f\wedge\partial_2f}{|\partial_1f\wedge\partial_2f|}\circ f^{-1} = \xi.
    \end{equation}
    Hence, the current defined in \eqref{eq:lem:orientability:current} is the one induced by the inclusion map $F\colon \Sigma\to\R^n$ of the orientable surface $\Sigma$.
\end{proof}

\section{Proof of the main results}\label{sec:proof}

In this section, we will provide the proofs of the affirmative statements in \Cref{thm:main_6pi,thm::main_8pi_curv,cor:KLS}. The examples implying sharpness are given in \Cref{sec:examples} below.

\begin{proof}[{Proof of \Cref{thm:main_6pi}}]
    By the monotonicity inequality \Cref{lem:ap:monotonicity}, the assumption $\mathcal W(\mu) < 6\pi$ implies $1\le\Theta^2(\mu,x)<3/2$ for all $x\in\spt\,\mu$. Moreover, $\spt\,\mu$ is connected by the Li--Yau inequality \eqref{eq:LY_closed}. Thus, \Cref{lem:density_6pi} yields $\Theta^2(\mu,x)=1$ for all $x\in \spt\,\mu$.
    Hence, the first part of the theorem follows from \Cref{lem:induced_surface}. Sharpness follows from \Cref{sec:double-bubble}.
\end{proof}

\begin{proof}[Proof of \Cref{thm::main_8pi_curv}]
    As before, the assumptions imply that for all $x\in\spt\,\mu$ we have $1\le\Theta^2(\mu,x)<2$, so $\Theta(\mu,x)=1$ by \Cref{lem:curv_var} and $\spt\,\mu$ is connected. The first part of the theorem thus follows from \Cref{lem:induced_surface} as above, noting that any embedded surface in $\R^3$ must be orientable.    
    Sharpness follows from \Cref{sec:branched_immersion}.
\end{proof}

\begin{proof}[Proof of \Cref{cor:KLS}]
    Again, the assumptions yield $1\le\Theta^2(\mu,x)<2$ for all $x\in\spt\,\mu$ and connectedness of $\spt\,\mu$. Thus, Lemmas \ref{lem:main_8pi} and \ref{lem:induced_surface} imply that $\mu$ is induced by a $W^{2,2}$-conformal embedding $F\in \mathcal{E}_\Sigma$ with $\Sigma$ orientable due to \Cref{,lem:orientability}. Applying the Gauss--Bonnet theorem for $W^{2,2}$-conformal immersions \cite[Remark 2.1]{MR2928715}, we have 
    \begin{align}
        \|B\|_{L^2(\mu)}^2 = 8\mathcal W(\mu) - 32\pi (1-p),
    \end{align}
    where $p\in \N_0$ is the genus of $\Sigma$. The control on $B$ is hence equivalent to controlling $p\leq C(\delta,n)$. Assume that $\mu_j$ is a sequence of varifolds as in the statement, induced by $F_j\in\mathcal{E}_{\Sigma_j}$ with $p_j=\operatorname{genus}\Sigma_j\to \infty$.
    By \cite{MR2928715,MR3276154}, we then have
    \begin{align}
        \liminf_{j\to\infty}\mathcal{W}(\mu_j) = \liminf_{j\to\infty}\mathcal{W}(F_j)\geq \liminf_{j\to\infty}\beta^n_{p_j} = 8\pi,
    \end{align}
    where we used \cite[Theorem 1.1]{KLS} in the last equality.
\end{proof}

\section{Applications}\label{sec:applications}

\subsection{Willmore minimization with prescribed boundary}

We now consider integral $2$-varifolds with a singular part of the first variation. This can be seen as a generalized notion of boundary by the divergence theorem
\begin{align}\label{eq:first_vari_bdry}
    \int \Div_\mu \Phi\, \mathrm{d}\mu  = - \int \langle\Phi, H\rangle \,\mathrm{d}\mu  + \int\langle \Phi,\nu\rangle\, \mathrm{d}\sigma,
\end{align}
for $\Phi\in C_c^1(\RR^n;\RR^n)$.
Here, $\nu=\nu_\mu$ is the generalized outer conormal and $\sigma=\sigma_\mu$ the generalized boundary measure.

Suppose that $\gamma \vcentcolon= \gamma^1\cup \dots\cup\gamma^\alpha\subset \R^3$ is the disjoint union of smooth embedded closed curves, $\alpha\in \N$. Let $\sigma_0 = \nu_0 m \mathcal{H}^1\mres \gamma$ be a vector valued Radon measure with $m\in L^\infty(\mathcal{H}^1\mres \gamma;\N)$ and $\nu_0\in L^\infty(\gamma;\mathbb{S}^2)$ such that $\nu_0(x) \perp T_x\gamma^k$ for $x\in \gamma^k$. Consider the variational problem
\begin{align}\label{eq:def_P}
	\mathcal{P}\vcentcolon = \inf\Big\{\mathcal{W}(\mu) \colon \mu\text{ satisfies }\eqref{eq:first_vari_bdry}\text{ with }\nu_\mu \sigma_\mu = \sigma_0, \;\spt\,\mu \text{ compact and connected}\Big\}.
\end{align}
In \cite[Theorem 4.1]{MR4141858}, the existence of a minimizer for \eqref{eq:def_P} was proven under the assumptions $\mathcal{P}<4\pi$ and $\alpha\geq 2$, however the same proof also yields existence in case $\alpha=1$. For a suitable class of boundary data, we have the following interior regularity result.

\begin{theorem}\label{thm:boundary_reg}
    Suppose that $\mathcal{P}<4\pi$ and
    \begin{align}\label{eq:LY_boundary_small}
        \mathcal{P}+2\sup_{x_0\in \RR^3\setminus \gamma} \int_{\gamma} \frac{\langle x-x_0, \nu_0(x)\rangle}{|x-x_0|^2}m(x)\,\mathrm{d}\mathcal{H}^1(x) <6\pi.
    \end{align}
    Then any minimizer for $\mathcal{P}$ 
    is real analytic in $\RR^3\setminus \gamma$.
\end{theorem}

It follows from the discussion at the end of \Cref{subsec:ex_circ_boundary} below that the set of admissible boundary data in \Cref{thm:boundary_reg} is nonempty.

\begin{remark}\label{rem:boundary}\leavevmode
    \begin{enumerate}        
        \item\label{item:tildeP} In view of \Cref{cor:KLS} and \Cref{lem:main_8pi}, we may also consider 
        \begin{align}
        	\tilde{\mathcal{P}} \vcentcolon = \inf \Big\{\mathcal{W}(\mu) \colon~& \mu\text{ satisfies }\eqref{eq:first_vari_bdry}\text{ with }\nu_\mu \sigma_\mu = \sigma_0, \;\spt\,\mu \text{ compact and connected, }\\ &  \mu\text{ is induced by an integral $2$-current $T$ in $\R^3\setminus \gamma$ with $\partial T=0$}\Big\}. 
        \end{align}
        Then, under the assumptions $\tilde{\mathcal{P}}<4\pi$ and
        \begin{align}
        	\tilde{\mathcal{P}}+2\sup_{x_0\in \RR^3\setminus \gamma} \int_{\gamma} \frac{\langle x-x_0, \nu_0(x)\rangle}{|x-x_0|^2}m\,\mathrm{d}\mathcal{H}^1(x) <8\pi,
        \end{align}
        we have existence of a minimizer in $\tilde{\mathcal{P}}$ which is real analytic in $\RR^3\setminus \gamma$.
        \item \Cref{thm:boundary_reg} remains valid if the assumption of connectedness of $\spt\,\mu$ in \eqref{eq:def_P} is removed or if $\RR^3$ is replaced by $\RR^n$, $n\geq 3$. Similarly for \Cref{rem:boundary}\eqref{item:tildeP}.
    \end{enumerate}
\end{remark}

\begin{proof}[Proof of \Cref{thm:boundary_reg}]
    By the monotonicity formula with boundary (see for instance \cite[p.~555]{MR4141858}), for any $\mu$ as in \eqref{eq:def_P} we have
    \begin{align}
        \Theta^2(\mu,x_0) \leq \frac{1}{4\pi}\mathcal{W}(\mu) + \frac{1}{2\pi} \int_\gamma \frac{\langle x-x_0,\nu_0(x)\rangle}{|x-x_0|^2}m(x)\,\mathrm{d}\mathcal{H}^1(x) \quad \text{ for all }x_0\in \RR^3\setminus \gamma,
    \end{align}
    so that the assumption \eqref{eq:LY_boundary_small} implies $\Theta^2(\mu,x_0)<3/2$. By \Cref{lem:density_6pi}, we conclude that $\Theta^2(\mu,x_0)=1$ for $x_0\in \spt\,\mu\setminus \gamma$. \Cref{thm:BiZhou_Scaled} yields that $\spt\,\mu\cap B_r(x_0)$ is parametrized by a conformal immersion $f\in W^{2,2}(U;\RR^3)$, $U\subset \R^2$ open. Since $\mu$ is minimizing for \eqref{eq:def_P}, we find that $\frac{\mathrm{d}}{\mathrm{d}t}\vert_{t=0}\mathcal{W}((f+t \Phi)(U))=0$ for all $\Phi\in C_c^\infty(U;\RR^3)$. Thus, $f$ is a solution to the weak Willmore equation in $U$ and hence real analytic by \cite[Theorem I.3]{MR2430975}.
\end{proof}

\subsection{Topology and Gauss--Bonnet for varifolds}\label{subsec:Gauss_Bonnet}

By the Gauss equations, we may express the Gauss curvature in terms of the mean curvature and the length of the second fundamental form. We use this to define the
\emph{Gauss curvature of an integral curvature varifold $\mu$} by
\begin{align}\label{eq:def_Gaussian}
    K \vcentcolon = \frac{|H|^2}{2} - \frac{|B|^2}{4}.
\end{align}

Motivated by the Gauss--Bonnet theorem, we find a notion of Euler characteristic for integral curvature varifolds.
\begin{definition}\label{def:chi}
    Let $\mu$ be an integral curvature varifold in $\R^n$. Suppose that $B\in L^2(\mu)$. Then we define the \emph{Euler characteristic of $\mu$} by
    \begin{align}
       \chi(\mu) \vcentcolon= \frac{1}{2\pi}\int K \,\mathrm{d}\mu.
    \end{align}
\end{definition}
It is not immediate whether \Cref{def:chi} gives an integer. In the setting of \Cref{thm:main_6pi,cor:KLS}, $\mu$ is already an integral curvature varifold with $B\in L^2(\mu)$ and we have the following.

\begin{cor}[Gauss--Bonnet theorem for varifolds]\label{cor:Gauss_Bonnet}
    Let $\mu$ be an integral $2$-varifold in $\R^n$ with $\mu(\RR^n)<\infty$.
    \begin{enumerate}
        \item  If $\mathcal{W}(\mu)<6\pi$, then $\chi(\mu)=\chi(\spt\,\mu)\in \Z_{\leq 2}$. 
        \item If $\mathcal{W}(\mu)<8\pi$ and $\mu$ is induced by an integral $2$-current $T$ in $\R^n$ with $\partial T=0$, then 
        $\chi(\mu) = 2-2p$ for $p\in\N_0$ where $p$ is the genus of $\spt\,\mu$.
    \end{enumerate}
\end{cor}
\begin{proof}
    This follows from \Cref{thm:main_6pi,cor:KLS}, the Gauss--Bonnet Theorem for $W^{2,2}$-immersions, see \cite[Remark 2.1]{MR2928715}, and by expressing the Euler characteristic in terms of the genus.
\end{proof}

We can now solve the Willmore problem in the class of varifolds.

\begin{proof}[Proof of \Cref{cor:equivalent_minimization}]
    By \cite{Kusner}, the infimum among smooth orientable surfaces of genus $p$ is strictly below $8\pi$.
    The equivalence of \ref{item:min_W_smooth} and \ref{item:min_W_conf} is known by \cite{Simon83,MR2928715,MR3276154}.

    For \ref{item:min_W_vari} and \ref{item:min_W_vari_current}, note that by \Cref{thm::main_8pi_curv,cor:KLS}, each varifold in the admissible set is induced by some immersion in $\mathcal{E}_\Sigma$. In case \ref{item:min_W_vari}, the Gauss--Bonnet Theorem (see \cite[Remark 2.1]{MR2928715}) implies that $\Sigma$ is orientable with genus $p$. For case \ref{item:min_W_vari_current}, this follows from \Cref{cor:KLS} and \Cref{cor:Gauss_Bonnet}. Thus \ref{item:min_W_vari} and \ref{item:min_W_vari_current} are equivalent to \ref{item:min_W_conf}.
\end{proof}

\subsection{Strong approximability}\label{subsec:approx}

It is a natural question whether (curvature) varifolds can be approximated by smooth surfaces with respect to the weak convergence of measures. While some results are known when considering subclasses given by graphs of sufficiently smooth functions, see \cite{MR1704565}, this is, in general, an open problem. In this section, we will see that approximation is possible even in a strong sense under suitable energy bounds.
The key observation is that a suitable mollification is well-behaved in the class of $W^{2,2}$-conformal Lipschitz immersions \cite[p.~316]{MR2928715}. However, the proof of $W^{2,2}$-convergence is much more technical than in Euclidean space due to the fact that the covariant differentiation $\nabla$ does not commute with the convolution. We thus decided to provide a full proof for completeness. 

We follow the approach in \cite[Section 4]{MR0710054} and use mollification in a tubular neighborhood. Since we want convergence in the best norm possible, we cannot simply apply  extension or trace operators.

\begin{lemma}\label{lem:mollification}
    Let $\Sigma$ be a closed smooth surface embedded in $\R^m$, $m\geq 3$. Let $\mathcal{U}\subset \R^m$ be a tubular neighborhood of $\Sigma$ and let $P\colon \mathcal{U}\to \Sigma$ be the nearest point projection. For any $F\in L^1(\Sigma;\R^N)$ we set $\bar{F}\vcentcolon = F\circ P$ and, for $\varepsilon>0$ small,
    \begin{align}
        \bar{F}_\varepsilon\colon \mathcal{U}_\varepsilon\to \R^N, \quad \bar{F}_\varepsilon(x) \vcentcolon= (\psi_\varepsilon*\bar{F})(x)\vcentcolon= \int_{\R^m} \psi_\varepsilon(x-y) \bar{F}(y)\,\mathrm{d}y.
    \end{align}
    Here $\mathcal{U}_\varepsilon \vcentcolon = \{ x\in \mathcal{U} \colon \dist(x,\partial\mathcal{U})>\varepsilon\}$ and $\psi_\varepsilon(x)\vcentcolon = \varepsilon^{-m}\psi(x/\varepsilon)$  for a smooth radial mollifier $\psi$ on $B_1(0)\subset \R^m$. Then 
    \begin{enumerate}[label=(\alph*)]
        \item\label{item:approx_1} $\bar{F}_\varepsilon \in C^\infty(\mathcal{U}_\varepsilon;\R^N)$;
        \item\label{item:approx_2}  if $F\in C(\Sigma;\R^N)$, then $\bar{F}_\varepsilon\to F$ uniformly on $\Sigma$;
        \item\label{item:approx_3}  if $F\in L^q(\Sigma;\R^N)$, $1\leq q<\infty$, then $\bar{F}_\varepsilon\to F$ in $L^q(\Sigma;\R^N)$;
        \item\label{item:approx_4}  if $F\in W^{1,q}(\Sigma;\R^N)$ then $\bar{F}\in W^{1,q}(\mathcal{U};\R^N)$ and componentwise we have
        \begin{align}\label{eq:DbarF}
             \mathrm{D}\bar{F}^j(y) = \mathrm{D}P(y)^T \nabla F^j(P(y))\qquad \text{for a.e.\ $y\in \mathcal{U}$, $j=1,\dots,N$.}
        \end{align}
    \end{enumerate}
\end{lemma}

\begin{proof}
    Without loss of generality, we may assume $N=1$.
    Items \ref{item:approx_1} and \ref{item:approx_2} follow directly from the properties of the mollification in $\R^m$. 
    Let $N\Sigma\vcentcolon=\{(p,v)\in \Sigma\times \R^m, v\perp T_p\Sigma\}\subset \R^m\times \R^m$ be the normal bundle of $\Sigma$. We assume that $\mathcal{U}$ is parametrized by the smooth diffeomorphism     
    \begin{equation}\label{eq:lem:approx:product_1}
        \varphi\colon N_\delta \vcentcolon = \{(p,v)\in N\Sigma \colon |v|<\delta\}\to\mathcal{U},\quad  \varphi(p,v) = p+v,
    \end{equation}
    with Jacobian $\theta(p,v)\vcentcolon=|\det \mathrm D\varphi(p,v)^*\circ \mathrm D\varphi(p,v)|^{1/2}$ which is bounded by positive numbers from above and below on $N_\delta$. Note that $P(\varphi(p,v))=p$ so the area formula yields
    \begin{align}\label{eq:Lp_vs_Lp}
        \int_{\mathcal{U}} |\bar{F}(y)|^q\,\mathrm{d}y = \int_{N_\delta} |F(p)|^q \theta(p,v)\,\mathrm{d}\mathcal{H}^m(p,v),
    \end{align}
    and thus $F\in L^q(\Sigma)$ if and only if $\bar{F}\in L^q(\mathcal{U})$. We also observe
    \begin{equation}
        \int_{\R^m}\psi_\varepsilon(x-y)\,\mathrm{d}y =1,\quad        
        \int_{\Sigma}\psi_\varepsilon(x-p)\,\mathrm d\mathcal{H}^2(p) \le C\varepsilon^{2-m} \qquad \text{for all $x\in\R^m$}.
    \end{equation}
    First, for $x\in \Sigma$, we note the pointwise estimate
    \begin{align}\label{eq:conv_ptwise}
        |\psi_\varepsilon*\bar{F}(x)| \leq \Bigl(\int_{\R^m} \psi_\varepsilon(x-y)\,\mathrm{d}y\Bigr)^{1-1/q} \Bigl( \int_{\Sigma(\varepsilon)} \psi_\varepsilon(x-y)|\bar{F}(y)|^q\,\mathrm{d}y\Bigr)^{1/q},
    \end{align}
    where $\Sigma(\varepsilon)\vcentcolon=\{y\in \mathcal{U} \colon \dist(y,\Sigma)<  \varepsilon\} = \varphi(N_\varepsilon)$ for $0<\varepsilon<\delta$,
    so that integration implies
    \begin{align}
        \int_\Sigma |\bar{F}_\varepsilon(x)|^q\mathrm{d}x &\leq \int_\Sigma \int_{\Sigma(\varepsilon)} \psi_\varepsilon(x-y)|\bar{F}(y)|^q\,\mathrm{d} y \,\mathrm{d}\mathcal{H}^2(x) \\
        &= \int_{N_\varepsilon} \int_\Sigma \psi_\varepsilon(x-\varphi(p,v))|F(p)|^q \theta(p,v) \,\mathrm{d}\mathcal{H}^2(x)\,\mathrm{d}\mathcal{H}^m(p,v) \leq C\Vert F\Vert_{L^q(\Sigma)}^q.\label{eq:conv_Lp_estimate}
    \end{align}
    Now, in order to prove \ref{item:approx_3}, we take $\gamma>0$ and $G\in C(\Sigma)$ with $\Vert F-G\Vert_{L^q(\Sigma)}<\gamma$. Using \eqref{eq:conv_Lp_estimate} applied to $(F-G)\circ P$ we estimate
    \begin{align}
        \Vert \bar{F}_\varepsilon - F\Vert_{L^q(\Sigma)} &\leq \Vert \psi_\varepsilon*(\bar{F}- \bar{G})\Vert_{L^q(\Sigma)} + \Vert \psi_\varepsilon* \bar{G} - G\Vert_{L^q(\Sigma)} + \Vert G-F\Vert_{L^q(\Sigma)} \\
        &\leq C\gamma + \Vert \psi_\varepsilon*\bar{G}-G\Vert_{L^q(\Sigma)} + \gamma.
    \end{align}    
    This can be made arbitrarily small, taking first $\gamma>0$ and then $\varepsilon>0$ small enough using \ref{item:approx_2}.
    For \ref{item:approx_4}, we take an approximating sequence $F_k\in C^\infty(\Sigma)$ with $\Vert F_k - F\Vert_{W^{1,q}(\Sigma)}\to 0$. Let $\bar{F_k}\vcentcolon= F_k\circ P\in C^\infty(\mathcal{U})$.    
    For any $\Phi\in C_c^\infty(\mathcal{U})$ and $i\in\{1,\dots,m\}$ we have
    \begin{align}
        -\int_{\mathcal{U}} \bar{F_k}\,\mathrm{D}_{e_i}\Phi\,\mathrm{d}y =    
        \int_{N_\delta} \langle \nabla F_k(p), 
        \mathrm{D}_{e_i} P(\varphi(p,v))\rangle \Phi(\varphi(p,v))\theta(p,v)\,\mathrm{d}\mathcal{H}^m(p,v).
    \end{align}
    Using $\bar{F_k}\to \bar{F}$ in $L^q(\mathcal{U})$ by \eqref{eq:Lp_vs_Lp} and $\nabla F_k \to \nabla F$ in $L^q(\Sigma;\R^m)$, it follows $\bar{F}\in W^{1,q}(\mathcal{U})$ with \eqref{eq:DbarF}. 
\end{proof}

The main advantage of this explicit global approximation procedure is that it preserves the property of being immersed. 

\begin{lemma}\label{lem:approx}
    Let $\Sigma$ be a closed smooth surface, and $F\in W^{2,2}(\Sigma;\RR^n)$ be a 
    Lipschitz immersion. Then the mollification $\bar F_\varepsilon \in C^\infty(\Sigma;\R^n)$ as in \Cref{lem:mollification} are immersions for $\varepsilon>0$ small and $\bar{F}_\varepsilon\to F$ strongly in $W^{2,2}(\Sigma;\RR^n)$ as $\varepsilon\to 0$.
\end{lemma}

\begin{proof}
    We may assume that $\Sigma$ is isometrically embedded in $\R^m$ for some integer $m\ge3$ and use the notation of \Cref{lem:mollification} and its proof in the sequel.  Let $T(p)\colon \R^m\to T_p\Sigma$ be the orthogonal projection, depending smoothly on $p\in \Sigma$. 
    For any $F\in C^\infty(\Sigma)$ with extension $\bar F \in C^\infty(\mathcal{U})$ the gradient and covariant Hessian (viewed as a vector in $\R^m$ or as a bilinear form on $T_p\Sigma\subset \R^m$, respectively) satisfy
    \begin{align}
        \nabla F (p)& = T(p) \mathrm{D}\bar{F}(p), \\
        \nabla^2 F(p)(\xi,\eta) 
        &= \langle T(p)\mathrm{D}_\xi T(p)\, \mathrm{D}\bar{F}(p),\eta\rangle + \langle T(p)\mathrm{D}_\xi\mathrm{D} \bar{F}(p),\eta\rangle,\label{eq:gradient_Hessian}
    \end{align}
    for $p\in\Sigma$, $\xi,\eta\in T_p\Sigma$.     
    For any $u\in W^{1,2}(\Sigma)$ and $p\in \Sigma$, by a Poincar\'e inequality, see \cite[(4.2)]{MR0710054}, we have
    \begin{align}
        \int_{B_\varepsilon^m(p)}|\bar u (x) -  \bar u_\varepsilon(p)|^2 \,\mathrm{d}x &\leq C\varepsilon^2 \int_{B_\varepsilon^m(p)}|\mathrm{D} \bar{u}(x)|^2\,\mathrm{d}x \\
        &\leq C \varepsilon^2 \int_{N_\varepsilon \cap (B_\varepsilon^m(p)\times \R^m)} |\mathrm{D}\bar{u}(\varphi(x,v))|^2 \theta(x,v)\,\mathrm{d}\mathcal{H}^m(x,v)\\
        &\leq C \varepsilon^m \int_{\Sigma\cap B_\varepsilon^m(p)}|\nabla u(x)|^2\,\mathrm{d}\mathcal{H}^2(x).\label{eq:poincare1}
    \end{align}
    Here we used the area formula and $|\mathrm{D}\bar{u}(y)| = |\nabla u (P(y))|$ for a.e.\ $y\in \mathcal{U}$ by \eqref{eq:DbarF}.     
    Taking $\nabla u$ componentwise, this estimate extends to $u\in W^{1,2}(\Sigma;\R^m)$. Since $F\in W^{2,2}(\Sigma;\R)$, by direct computation,  $u\vcentcolon=\nabla F\in W^{1,2}(\Sigma;\R^m)$ with
    \begin{align}
        \nabla_\xi \langle u, e_j\rangle = \nabla^2 F(\xi,T e_j)+ \langle \nabla F, \mathrm{D}_\xi T\;e_j\rangle.
    \end{align}
    We conclude that for $F\in \mathcal{E}_\Sigma$ we may take $u(p)\vcentcolon = (p,\nabla F(p)) \colon \Sigma \to \R^m\times  \R^{n\times m}$ in \eqref{eq:poincare1} and find
    \begin{align}\label{eq:lem:approx:poincare_0}
        \int_{B_\varepsilon^m(p)}|\bar{u}(x) -  \bar u_\varepsilon(p)|^2 \,\mathrm{d}x &\leq C\varepsilon^m \int_{\Sigma\cap B_\varepsilon^m(p)} (1+|\nabla F(x)|^2+|\nabla^2 F(x)|^2)\,\mathrm{d}\mathcal{H}^2(x).
    \end{align}
   For $\lambda>0$, define $M_\lambda\subset \Sigma\times \R^{n\times m}$ as the open set of all tuples $(p,L)$ such that $p\in \Sigma$, $L\in \R^{n\times m}$, and $\det ((L\vert_{T_p\Sigma})^*\circ L\vert_{T_p\Sigma})>\lambda$. Since $F\in \mathcal{E}_\Sigma$, there exists $\lambda>0$ such that $u(x)\in M_\lambda$ for $\mathcal{H}^2$-a.e.\ $x\in \Sigma$, and consequently, by the area formula with $\varphi$ as in \eqref{eq:lem:approx:product_1}, $\bar{u}(x)\in M_\lambda$ for $\mathcal{L}^m$-a.e.\ $x\in \mathcal{U}$. It follows from \eqref{eq:lem:approx:poincare_0}
    that for all $p\in \Sigma$ we have
    \begin{equation}\label{eq:lem:approx:poincare_1}
                \dist(\bar u_\varepsilon(p),M_\lambda)^2\le C\varepsilon^{-m}\int_{B^m_\varepsilon(p)}|\bar u(x) - \bar u_\varepsilon(p)|^2\,\mathrm dx \leq C \int_{\Sigma\cap B_\varepsilon^m(p)} (1+|\nabla F|^2+|\nabla^2 F|^2)\,\mathrm{d}\mathcal{H}^2.
    \end{equation}
    For $\varepsilon>0$ small enough, we thus find $\bar{u}_\varepsilon(p)\in M_{\lambda/2}$. 
    Using \eqref{lem:mollification}\ref{item:approx_4} and the definition of the gradient, we note $\nabla \bar{F}_\varepsilon(p) = T(p) (\psi_\varepsilon*\mathrm{D}\bar{F})(p)$, $p\in\Sigma$.
    We further note that $\mathrm{D}P(p)=T(p)$ for $p\in \Sigma$ and $T(p)^T = T(p)$. Thus, by \eqref{eq:DbarF}, for any $y\in \mathcal{U}$ we have
    \begin{align}
        &T(p) \mathrm{D}\bar{F}(y) - \nabla F(P(y)) = T(p) \mathrm{D}P(y)^T \nabla F(P(y))- T(P(y)) \mathrm{D}P(P(y))^T \nabla F(P(y)) \\
        &\qquad = [T(p)-T(P(y))] \mathrm{D}P(y)^T \nabla F(P(y)) + T(P(y)) [\mathrm{D}P(y)^T-\mathrm{D}P(P(y))^T]\nabla F(P(y)). \label{eq:approx_error_nablaF_0}
    \end{align}
    Using the Lipschitz continuity of $T$ and $\mathrm{D}P$ we find
    \begin{align}
        |T(p)\mathrm{D}\bar{F}(y)-(\nabla F)(P(y))| \leq C |\nabla F(P(y))| (|p-P(y)| +|y-P(y)|).\label{eq:approx_error_nablaF_1}
    \end{align}
   By \eqref{eq:conv_ptwise}, we thus have
    \begin{align}
        |\nabla \bar{F}_\varepsilon(p)- \pi_2 \bar{u}_\varepsilon(p)| &=
       |T(p)\,(\psi_\varepsilon*\mathrm{D}\bar{F})(p) - (\psi_\varepsilon* (\nabla F\circ P))(p)|\\
       &\leq \int_{B_\varepsilon^m(p)}\psi_\varepsilon(p-y) |T(p)\mathrm{D}\bar{F}(y)-(\nabla F)(P(y))|\,\mathrm{d}y \leq C\varepsilon \int_\mathcal{U}|\nabla F(P(y))|\,\mathrm{d}y\label{eq:approx_error_nablaF_2}
    \end{align}
    where we used $|y-P(y)|\leq \varepsilon$, $|p-P(y)|\leq \varepsilon$ in the last step. This goes to zero as $\varepsilon\to 0$ by \eqref{eq:Lp_vs_Lp}, and thus $\bar F_\varepsilon|_{\Sigma}$ is a smooth immersion for $\varepsilon>0$ small enough.
    
    It remains to show that $\bar{F}_\varepsilon\to F$ in $W^{2,2}(\Sigma)$. First, $\bar{F}_\varepsilon\to F$ in $L^2(\Sigma)$ follows directly from \Cref{lem:mollification}\ref{item:approx_3}. For the second derivative, we consider the components $F^j$ of $F$, $j=1,\dots,n$. 
            By \Cref{lem:mollification}\ref{item:approx_4}, $\bar{F}^j\in W^{2,2}(\Sigma)$ and for $p\in \Sigma$ and $\eta,\xi\in T_p\Sigma$, $|\eta|,|\xi|\leq 1$, we compute
            \begin{align}
                \nabla^2 \bar{F}^j_\varepsilon(p) (\xi,\eta) = 
                \int_{B_\varepsilon^m(p)} \psi_\varepsilon(p-y)&\Big[\langle T(p) \mathrm{D}_\xi T(p) \mathrm{D}\bar{F}^j(y),\eta\rangle +\langle T(p) \mathrm{D}_\xi\mathrm{D}\bar{F}^j(y),\eta\rangle\Big]\,\mathrm{d}y.\label{eq:approx_nabla^2_1}
            \end{align}
            We compare this with $\psi_\varepsilon*(\nabla^2 F^j\circ P)$ since by \Cref{lem:mollification}\ref{item:approx_3} we have $\psi_\varepsilon*(\nabla^2 F^j\circ P)\to \nabla^2 F$ in $L^2(\Sigma;\R^{m\times m})$. We have
            \begin{align}
                (\psi_\varepsilon*(\nabla^2F^j\circ P))(p)(\xi,\eta) = \int_{B_\varepsilon^m(p)} \psi_\varepsilon(p-y) &\Big[\langle T(P(y)) \mathrm{D}_\xi T(P(y))\mathrm{D} \bar{F}^j(P(y)),\eta\rangle \\
                &\qquad\qquad + \langle T(P(y))\mathrm{D}_\xi \mathrm{D}\bar{F}^j(P(y)),\eta\rangle\Big]\,\mathrm{d}y.\label{eq:approx_nabla^2_2}
            \end{align}
            To compare this to \eqref{eq:approx_nabla^2_1}, we differentiate \eqref{eq:DbarF} to find
            \begin{align}
                \mathrm{D}_\xi\mathrm{D} \bar{F}^j(y) &= \mathrm{D}_\xi \mathrm{D}P(y)^T T(P(y)) \mathrm{D}\bar F^j (P(y)) + \mathrm{D}P(y)^T \mathrm{D}T(P(y)) \mathrm{D}_\xi P(y) \mathrm{D}\bar{F}^j(P(y)) \\
                &\qquad + \mathrm{D}P(y)^T T(P(y)) \mathrm{D}^2\bar{F}^j(P(y)) \mathrm{D}_\xi P(y). \label{eq:D^2barF}
            \end{align}
            Applying this to the highest order term in \eqref{eq:approx_nabla^2_1} and using $P(P(y))=P(y)$, we obtain
            \begin{align}
               \mathrm{D}_\xi \mathrm{D}\bar{F}^j(P(y)) & =\mathrm{D}_\xi \mathrm{D}P(P(y))^T T(P(y)) \mathrm{D}\bar F^j (P(y)) + \mathrm{D}P(P(y))^T \mathrm{D}T(P(y)) \mathrm{D}_\xi P(P(y)) \mathrm{D}\bar{F}^j(P(y)) \\
                &\qquad + \mathrm{D}P(P(y))^T T(P(y)) \mathrm{D}^2\bar{F}^j(P(y)) \mathrm{D}_\xi P(P(y)).\label{eq:D^2barF}
            \end{align}
            Therefore, adding and subtracting as in \eqref{eq:approx_error_nablaF_0} and using the Lipschitz properties of $T,P$ and their derivatives as in \eqref{eq:approx_error_nablaF_1},  we find
            \begin{align}
                &|\langle T(p) \mathrm{D}_\xi \mathrm{D}\bar{F}^j(y),\eta\rangle- \langle T(P(y)) \mathrm{D}_\xi \mathrm{D}\bar{F}^j(P(y)), \eta\rangle | \\
                &\qquad \leq C \Big( |\mathrm{D}\bar F^j(P(y))|+|\mathrm{D}^2 \bar{F}^j(P(y))|\Big)\big(|y-P(y)|+|p-P(y)|\big)
            \end{align}
            recalling that $|\eta|,|\xi|\leq 1$. Proceeding similarly for the first order terms in \eqref{eq:approx_nabla^2_1} and \eqref{eq:approx_nabla^2_2}, in analogy to \eqref{eq:conv_Lp_estimate} we have
            \begin{align}
                &\Vert \nabla^2 \bar{F}^j_\varepsilon - \psi_\varepsilon*(\nabla^2F^j\circ P)\Vert_{L^2(\Sigma)}^2 \leq C \varepsilon^2\int_{\mathcal{U}} (|\mathrm{D}\bar F^j(P(y))|^2 + |\mathrm{D}^2\bar F^j(P(y))|^2) \,\mathrm{d}y.\label{eq:approx_nabla^2_3}
            \end{align}
            Recall from \eqref{eq:DbarF} that $|\mathrm{D}\bar{F}^j(P(y))| = |\nabla F^j(P(y))| \in L^2(\mathcal{U})$ by \eqref{eq:Lp_vs_Lp}. Therefore, to show integrability of the second order term in \eqref{eq:approx_nabla^2_3} it suffices to consider the second order term in $\mathrm{D}^2 \bar F^j(P(y))(\xi,\eta)$, as given by \eqref{eq:D^2barF}, $y\in \mathcal{U}$ and $\xi,\eta\in \R^m$, which is
            \begin{align}
                \langle \mathrm{D}P(P(y))^T T(P(y)) \mathrm{D}^2 \bar{F}^j(P(y)) \mathrm{D}_\xi P(P(y)),\eta\rangle = \langle \mathrm{D}^2\bar{F}^j(P(y)) T(P(y))\xi, T(P(y))\eta\rangle.
            \end{align}
            Hence, this term is controlled by the norm of $\mathrm{D}^2\bar{F}^j(P(y))\vert_{T_{P(y)}\Sigma\times T_{P(y)}\Sigma}$. This is the second order term in the Hessian, cf.\ \eqref{eq:gradient_Hessian}, which leads to the estimate
            \begin{align}
                |\mathrm{D}^2\bar F^j(P(y))| \leq C\left(|\nabla F^j(P(y))| + |\nabla^2 F^j(P(y))|\right).
            \end{align}
            Consequently, \eqref{eq:approx_nabla^2_3}, \eqref{eq:Lp_vs_Lp}, and $F^j \in W^{2,2}(\Sigma)$ imply $\Vert \nabla^2 \bar{F}_\varepsilon^j - \psi_\varepsilon*(\nabla^2 F^j \circ P)\Vert_{L^2(\Sigma)} \to 0$. This completes the proof of $\bar F_\varepsilon \to F$ in $W^{2,2}(\Sigma;\R^n)$.
\end{proof}

\begin{cor}[Strong approximability]\label{cor:approximation}
    Let $\mu$ be an integral curvature $2$-varifold in $\R^n$ with $\mu(\R^n)<\infty$, $B\in L^2(\mu)$, and $\mathcal{W}(\mu)<8\pi$. Then there exist a compact smooth surface $\Sigma$ without boundary with $\Sigma\cong \spt\,\mu$ and a sequence $F_j\in C^\infty(\Sigma;\RR^n)$ of embeddings approximating $\mu$ in the $W^{2,2}(\Sigma;\RR^n)$-norm. 
\end{cor}
\begin{proof}
    By \Cref{thm::main_8pi_curv}, $\mu$ is induced by some $F\in \mathcal{E}_\Sigma$ for some compact smooth surface $\Sigma$ without boundary. \Cref{lem:approx} gives a sequence of smooth immersions $F_j$ strongly approximating $F$, and the Li--Yau inequality implies that $F_j$ is an embedding for $j\in\N$ large.
\end{proof}

Let $\Omega\subset \R^n$ be an open set. We recall from \cite{MR1704565} the classes $F_2 C^\infty(\Omega)$ and $F_2 W^{2,p}(\Omega)$ of varifolds in $\Omega$ consisting of locally finite unions of $2$-dimensional graphs of functions of the class $C^\infty$ and $W^{2,p}$, $p>2$, respectively. 
The class of $2$-dimensional integral curvature varifolds in $\Omega$ with $B\in L^p$ is denoted by $\mathcal{C}_2^p(\Omega)$. For $A>0$ and $\Lambda>0$, consider the set
\begin{align}
    E \vcentcolon = \{ \mu \in \mathcal{C}_2^p(\Omega) \colon \mu(\Omega) \leq A, \Vert B\Vert_{L^p(\mu)}\leq \Lambda \}.
\end{align}
For $p>2$, it was proven in \cite{MR1704565} that, for the weak-$*$-closure in the space of Radon measures, we have
\begin{align}
    \overline{E\cap F_2 C^\infty(\Omega)}\not\subset E\cap F_2W^{2,p}(\Omega) \subsetneq E.
\end{align}
To put \Cref{cor:approximation} in context with this result, denote by $I_2 C^\infty(\R^n)$ and $I_2W^{2,2}(\R^n)$ the class of integral curvature varifolds in ${\mathcal C}_2^2(\R^n)$ consisting of locally finite unions of  immersed $2$-manifolds of class $C^\infty$ and $W^{2,2}$, respectively. For $A>0$, $\Lambda$, $\delta>0$, we set
\begin{align}
    \tilde{E} \vcentcolon = \{ \mu \in {\mathcal{C}}_2^2(\R^n) \colon \mu(\R^n) \leq A, \Vert B\Vert_{L^2(\mu)}\leq \Lambda, \mathcal{W}(\mu)\leq 8\pi-\delta\}.
\end{align}
Then, \Cref{cor:approximation} and \Cref{thm::main_8pi_curv} imply
\begin{align}
   \overline{\tilde{E}\cap I_2C^\infty(\R^n)}  = \tilde{E}\cap I_2W^{2,2}(\R^n)= \tilde{E}.
\end{align}

\subsection{Adaptation to the Helfrich energy}\label{subsec:Helfrich}

We recall the definition of the Helfrich energy for oriented varifolds, see \cite[Section 2.2]{ChrisFabian}. An \emph{oriented $2$-varifold in $\RR^3$} is a Radon measure $V$ over $\mathbb{G}_2^{\mathrm{o}}(\RR^3)\vcentcolon= \RR^3 \times \mathbb G^{\mathrm{o}}(3,2)$, where $\mathbb G^{\mathrm{o}}(3,2)$ is the set of oriented $2$-dimensional subspaces in $\RR^3$ \cite{MR0825628}. We say that $V$ is \emph{integral} if 
\begin{align}
    V(k) = \int_M \Big(k(x,\xi(x)) \theta_1(x)+k(x,-\xi(x))\theta_2(x)\Big)\,\mathrm{d}\mathcal{H}^2(x) \qquad \text{for all }k\in C_c^0(\mathbb{G}_2^{\mathrm{o}}(\RR^3)).
\end{align}
Here $M\subset \R^3$ is an $\mathcal H^2$-rectifiable set, $\theta_1,\theta_{2}\in L^1_{\mathrm{loc}}(\mathcal{H}^2\mres M; \N)$, and $\xi$ is measurable, $|\xi(x)|=1$, and $\xi(x)$ orients $T_x\mu_V$ for $\mathcal{H}^2\mres M$-a.e.\ $x\in \R^3$ where $\mu_V\vcentcolon=(\theta_1+\theta_2) \mathcal{H}^2\mres M$ is the \emph{weight measure}. Any Lipschitz immersion $F\in \mathcal{E}_{\Sigma}$ of an oriented surface induces an integral oriented $2$-varifold, see \cite[Example 2.4]{ChrisFabian}. 
With $H=H_{\mu_V}$ as in \eqref{eq:first_vari}, we define the \emph{Helfrich energy} as
\begin{align}
    \mathcal{H}_{c_0}(V) \vcentcolon= \frac{1}{4} \int_{\mathbb{G}_2^{\mathrm{o}}(\RR^3)}|H(x)-c_0(\star\xi)|^2\,\mathrm{d} V(x,\xi).
\end{align}
Here the parameter $c_0\in\R$ is called \emph{spontaneous curvature} and $\star$ denotes the Hodge star operator. 
We also define the \emph{concentrated volume of $V$ at $x_0\in\RR^3$} by
\begin{align}
    \mathcal{V}_c(V,x_0) \vcentcolon = - \int_{\mathbb{G}_2^{\mathrm{o}}(\RR^3)} \frac{\langle x-x_0,\star\xi\rangle}{|x-x_0|^2}\,\mathrm{d}V(x,\xi),
\end{align}
see \cite[Section 3]{ChrisFabian}. 
Analogously to the Willmore functional, also the Helfrich energy satisfies a Li--Yau-type inequality \cite{ChrisFabian}.
Since the Helfrich functional naturally incorporates orientation, we impose a similar assumption as in \Cref{cor:KLS}. 

\begin{theorem}\label{thm:helfrich_8pi}
    Let $c_0\in \RR$.  Let $V$ be an integral oriented $2$-varifold in $\R^3$ with $\mu_V(\RR^3)+\mathcal{H}_{c_0}(V)<\infty$. Suppose that
    \begin{align}\label{eq:thm:helfrich_8pi}
        \mathcal{H}_{c_0}(V) + 2 c_0 \mathcal{V}_c(V,x_0) < 8\pi \qquad & \text{for all }x_0\in\RR^3,\\
    \label{eq:dT=0}
        \int_{\mathbb{G}_2^{\mathrm{o}}(\RR^3)} \langle \operatorname{curl} \Phi(x), \star\xi\rangle\, \mathrm{d}V(x,\xi) = 0 \qquad 
        & \text{for all } \Phi\in C_c^\infty(\RR^3;\RR^3).
    \end{align}
     Then $V$ and $\mu_V$ are induced by a conformal embedding $F\in \mathcal{E}_\Sigma$ of a compact Riemann surface $(\Sigma, g_0)$ without boundary.
\end{theorem}

\begin{remark}\leavevmode
    \begin{enumerate}
        \item Suppose $V$ is a \emph{varifold with enclosed volume} in the sense of \cite{ChrisFabian}. Then \cite[(4.19)]{ChrisFabian} directly implies \eqref{eq:dT=0} since $\Div \circ \operatorname{curl}=0$.
        \item There exists a universal constant $C$ such that for $c_0<0$, \eqref{eq:thm:helfrich_8pi} is implied by the sufficient condition
        \begin{equation}
            \mathcal H_{c_0}(V)\le 4\pi\Bigl(1 + \sqrt{1 + \frac{v_0|c_0|}{a_0C}}\Bigr) 
        \end{equation}
        where $a_0\vcentcolon=\mu_V(\R^3)$ and $v_0\vcentcolon=-\frac{1}3\int_{\mathbb G_2^\mathrm{o}(\R^3)}\langle x,\star\xi\rangle\,\mathrm dV(x,\xi)$, see \cite[Theorem 5.3]{scharrer2023properties}. Similarly, for $c_0\ge0$, \eqref{eq:thm:helfrich_8pi} is implied by the sufficient condition
        \begin{equation}
            \mathcal H_{c_0}(V)\le 8\pi -6c_0(4\pi^2 v_0)^\frac{1}{3},
        \end{equation}
        see \cite[Remark 6.11(ii)]{ChrisFabian}
    \end{enumerate}
\end{remark}

\begin{proof}[{Proof of \Cref{thm:helfrich_8pi}}]
    Let $T$ be the $2$-current given by
    \begin{align}
        T(\omega) \vcentcolon= \int_{M} \langle \omega(x),\xi(x)\rangle (\theta_1(x)-\theta_2(x))\,\mathrm{d}\mathcal{H}^2(x) = \int_{\mathbb{G}_2^{\mathrm{o}}(\RR^3)} \langle \omega(x), \xi\rangle \,\mathrm{d}V(x,\xi)
    \end{align}
    for two forms $\omega$ on $\RR^3$. Recalling that for a one form $\omega(x) = \sum_{k=1}^3\Phi_k(x) \mathrm{d}x_k$, $x\in\RR^3$, we have $\star\mathrm{d}\omega(x) = \operatorname{curl}\Phi(x)$, condition \eqref{eq:dT=0} is clearly equivalent to $\partial T=0$. 
    The Li--Yau inequality for the Helfrich energy \cite[Corollary 4.11]{ChrisFabian} implies
    \begin{align}
        \Theta^2(\mu_V,x_0) \leq \frac{1}{4\pi}\mathcal{H}_{c_0}(V) + \frac{c_0}{2\pi}\mathcal{V}_c(V,x_0) <2.
    \end{align}
    Since we also have $\mathcal{W}(\mu)<\infty$, we conclude $\Theta^2(\mu_V,x)=1$ for all $x\in\spt\,\mu_V$ by \Cref{lem:main_8pi}. Without loss of generality, we may thus assume $\theta_1=1$ and $\theta_2=0$ $\mu_V$-a.e.    
    From \Cref{lem:orientability}, it follows that $\mu_V$ and $T$ are induced by a conformal embedding $F$ as desired where $\xi$ is given by \eqref{eq:lem:orientability:normal}. It follows that also $V$ is induced by $F$ in the sense of \cite[Example 2.4]{ChrisFabian}.
\end{proof}



\section{Examples}\label{sec:examples}

\subsection{The standard double bubble}\label{sec:double-bubble}

\begin{figure}
     \centering
     \hfill
     \begin{subfigure}[b]{0.45\textwidth}
         \centering
        \includegraphics[height=5cm]{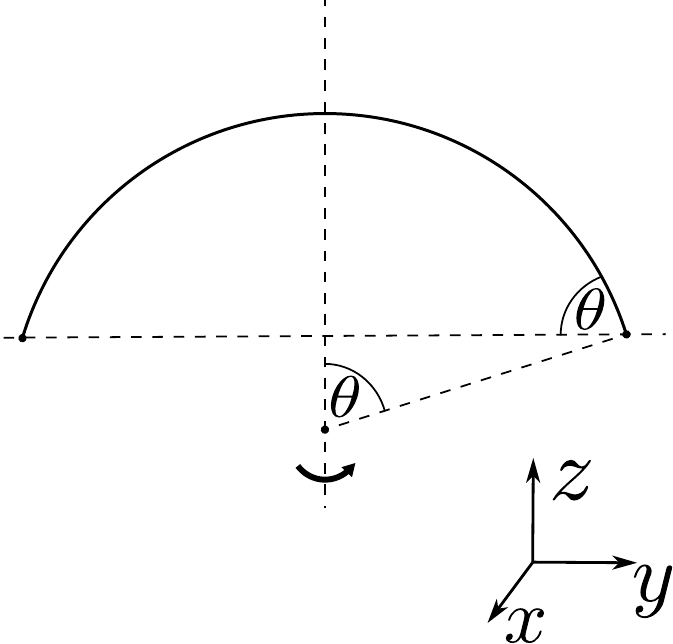}
         \caption{A rotationally symmetric spherical cap.}
          \label{fig:cap}
     \end{subfigure}
     \hfill
     \begin{subfigure}[b]{0.45\textwidth}
         \centering
         \includegraphics[height=5cm]{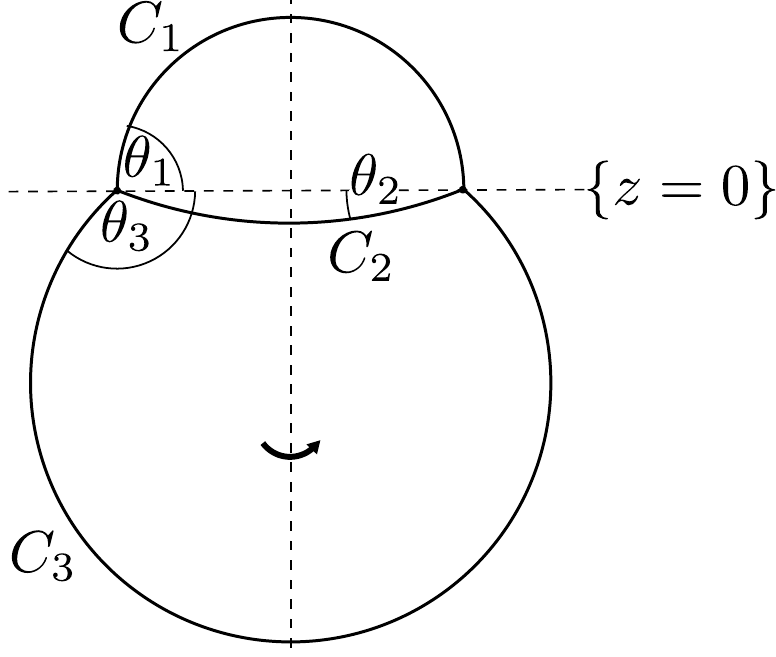}
         \caption{A standard double bubble.}
         \label{fig:double_bubble}
     \end{subfigure}
     \caption{Computing the Willmore energy of a double bubble.}
\end{figure}

Let $C \subset \RR^3$ be a closed spherical cap with radius $R>0$ and opening angle $\theta \in [0, \pi]$ (such that for $\theta =\pi/2$ we obtain the hemisphere). Assume that $C$ is obtained by revolving an axially symmetric arc in the plane $\{x=0\}$ around the $z$-axis. Elementary trigonometry implies that the intersection angle between the spherical cap and its base disk is also given by $\theta$, see \Cref{fig:cap}.

A famous configuration from minimal surface theory is a \emph{standard double bubble} \cite{MR1906593}, consisting of three spherical caps $C_1,C_2,C_3$, intersecting at 120\textdegree~along common boundary circle lying in the plane $\{z=0\}$, see \Cref{fig:double_bubble}. 
Thus, the opening angles of the spheres satisfy the conditions
\begin{align}\label{eq:angles}
    \theta_1 + \theta_2 = \frac{2\pi}{3},\qquad \theta_3-\theta_2 = \frac{2\pi}{3},\qquad \theta_1 + \theta_3 = \frac{4\pi}{3}.
\end{align}
Moreover, the radii are determined by the relation $\frac{1}{R_1} = \frac{1}{R_2}+ \frac{1}{R_3}$. By standard formulas, the area of the spherical caps is given by
\begin{align}
    \mathcal{H}^2(C_i) = 2\pi R_i^2 (1-\cos\theta_i).
\end{align}
Let $\Sigma \vcentcolon = C_1\cup C_2\cup C_3$. Then,  $\mu\vcentcolon=\mathcal{H}^2\mres \Sigma$ has locally bounded first variation, and the boundary terms cancel, \cite[Lemma 3.1]{MR1906593}, so that $\mu$ satisfies \eqref{eq:first_vari}. The Willmore energy is given by
\begin{align}
    \mathcal{W}(\mu) = \frac{1}{4} \sum_{i} \Big(\frac{2}{R_i}\Big)^2 \mathcal{H}^2(C_i) = \sum_i 2\pi(1-\cos\theta_i) = 6\pi - \sum_i\cos\theta_i = 6\pi - \cos\theta_2\big(1+2\cos \frac{2\pi}{3}\big),
\end{align}
where we used \eqref{eq:angles} in the last step. Since $\cos\frac{2\pi}{3} = -\frac{1}{2}$, we find that $\mathcal{W}(\mu)=6\pi$. Along the circle $\Sigma\cap \{z=0\}$, the blow up is given by the union of three half spaces, corresponding to the revolution of three arcs joining, and thus we have $\Theta^2(\mu,\cdot)=3/2$.

\subsection{The standard triple bubble}\label{sec:triple-bubble}
\begin{figure}
     \centering
     \hfill
     \begin{subfigure}[t]{0.4\textwidth}
         \centering
         \includegraphics[height=5cm]{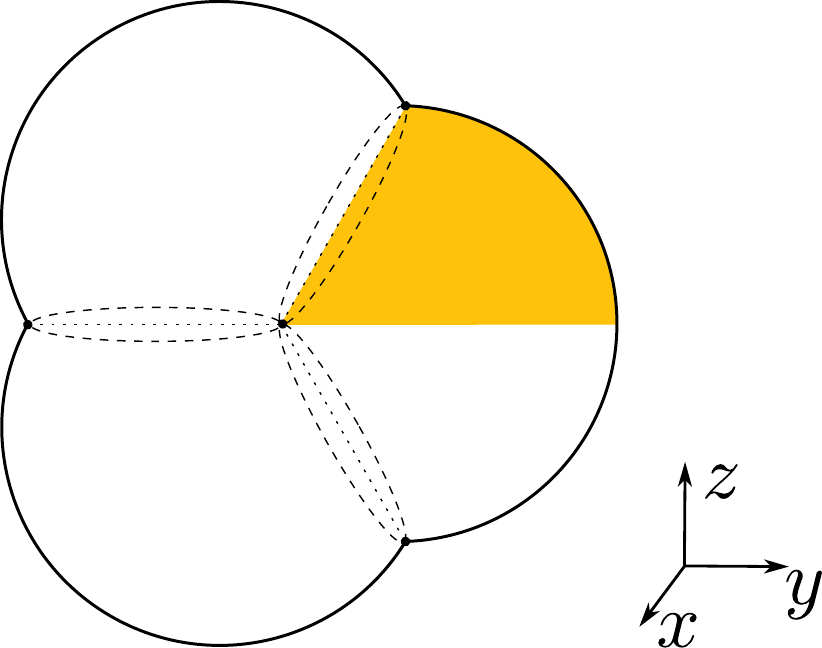}
    \caption{Profile of a symmetric triple bubble through the plane of symmetry $\{x=0\}$.}
          \label{fig:triple_bubble_1}
     \end{subfigure}%
     \hspace{0.1\textwidth}
      \begin{subfigure}[t]{0.4\textwidth}
         \centering
         \includegraphics[height=5cm]{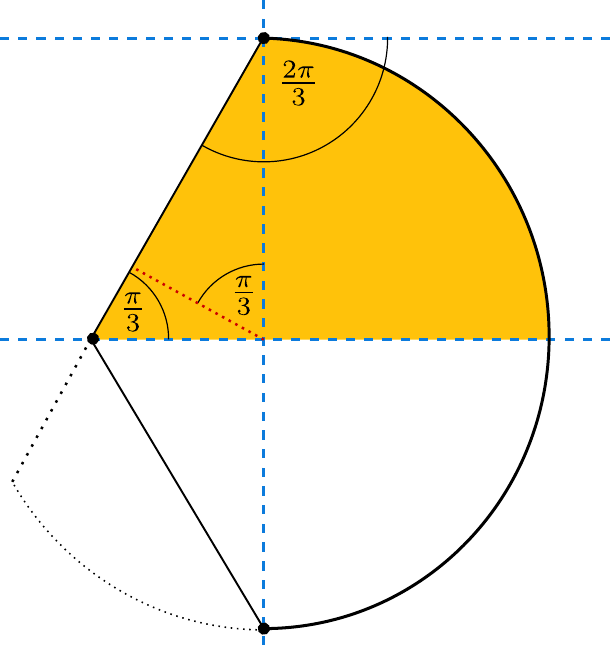}
         \caption{Determination of the opening angles.}
         \label{fig:triple_bubble_2}
     \end{subfigure}
     \hfill
     \\
     \hfill
     \begin{subfigure}[b]{0.9\textwidth}
         \centering
         \includegraphics[height=5cm,trim= {2cm 3.5cm 2cm 2cm},clip]{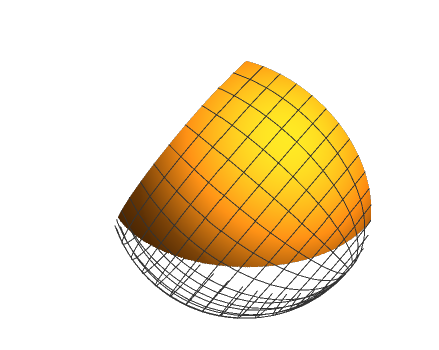}
         \caption{$3$-D plot of the spherical part parametrized by $X$. The grid corresponds to $\theta\in (\pi/3,\pi)$, $\varphi\in (-\pi,0)$, whereas the colored part corresponds to $(\theta,\varphi)$ as in \eqref{eq:theta_phi_range}.}
         \label{fig:triple_bubble_3}
     \end{subfigure}
     \hfill
    \caption{Parametrization of the symmetric triple bubble.}
\end{figure}

We consider a symmetric \emph{triple bubble} \cite{milman2024structure} consisting of three congruent parts of the unit sphere and three congruent parts of a disk, cf.\ \Cref{fig:triple_bubble_1}. We are going to compute its Willmore energy which coincides with the area of the three spherical parts. Each spherical part is symmetric with respect to two planes. Hence, it will be enough to compute the area of just one quarter of one spherical part. Consider a sphere with a removed spherical cap of opening angle $\alpha = \pi/3$:

\begin{equation}
    \tilde{\tilde{X}}(\theta,\varphi) \vcentcolon= 
    \begin{pmatrix}
        \sin \theta \cos\varphi \\
        \sin \theta \sin\varphi \\
        \cos\theta
    \end{pmatrix}
    \qquad \theta\in (\frac{\pi}{3},\pi), \quad \varphi\in (-\pi,\pi).
\end{equation}
In order to have the two symmetry planes given by $\{x=0\}$ and $\{z=0\}$, we first rotate $\tilde{\tilde{X}}$ by $\pi/2$ around the $z$-axis:
\begin{equation}
    \tilde{X}(\theta,\varphi) \vcentcolon= 
    \begin{pmatrix}
        -\sin \theta \sin\varphi \\
        \sin \theta \cos\varphi \\
        \cos\theta
    \end{pmatrix}
\end{equation}
and then rotate $\tilde{X}$ by 
$\pi/3$ around the $x$-axis as in \Cref{fig:triple_bubble_2}
\begin{equation}
    X(\theta,\varphi)\vcentcolon=
    \begin{pmatrix}
        -\sin\theta\sin\varphi \\
        \frac{1}2\sin\theta\cos\varphi - \frac{\sqrt{3}}2\cos\theta \\
        \frac{\sqrt{3}}2\sin\theta\cos\varphi + \frac{1}2\cos\theta
    \end{pmatrix}.
\end{equation}

In order to obtain just one quarter of the spherical part, we need the first and the third component of $X$ to be positive, cf.\ \Cref{fig:triple_bubble_3}. We thus arrive at the conditions 
\begin{equation}
    \frac{\sqrt{3}}2\sin\theta\cos\varphi + \frac{1}2\cos\theta>0,\qquad \varphi \in (-\pi,0)
\end{equation}
which is only nonempty for $\theta\in(\frac{\pi}3,\pi)$ if $\theta\in(\frac{\pi}3,\frac{5\pi}{6})$. Summarizing, we have
\begin{equation}\label{eq:theta_phi_range}
    \theta\in(\frac{\pi}3,\frac{5\pi}{6}),\qquad\varphi \in (-\varphi_\theta,0),\quad \varphi_\theta\vcentcolon=\arccos\Bigl(-\frac{1}{\sqrt{3}}\cot\theta\Bigr).
\end{equation}

The induced metric tensor $G$ satisfies $\sqrt{\det G}(\theta,\varphi) = \sin\theta.$
Let $\mu$ be the varifold associated with the triple bubble. Using the substitution $\theta = \arccos x$, it follows
\begin{align}
    \frac{\mathcal  W(\mu)}{12}&=\int_{\frac{\pi}3}^{\frac{5\pi}6}\int_0^{\arccos(-\frac{1}{\sqrt{3}}\cot(\theta))}\sqrt{\det G}(\theta,\varphi)\,\mathrm d\varphi\,\mathrm d\theta\\
    &=\int_{\frac{\pi}3}^{\frac{5\pi}6}\sin(\theta)\arccos\Bigl(-\frac{1}{\sqrt{3}}\cot\theta\Bigr)\,\mathrm d\theta = \int_{-\frac{\sqrt{3}}2}^{\frac{1}2}\arccos\Bigl(-\frac{x}{\sqrt{3-3x^2}}\Bigr) \mathrm dx.
\end{align}
Using the substitution $y=-\frac{x}{\sqrt{3-3x^2}}$ and subsequent integration by parts, we compute further
\begin{align}
    \frac{\mathcal W(\mu)}{12} & = \int_{-\frac{1}3}^1\frac{\sqrt{3}\arccos y}{(1+3y^2)^\frac{3}{2}}\,\mathrm dy  = \sqrt{3}\int_{-\frac{1}3}^1\frac{y}{\sqrt{1+3y^2}}\frac{\mathrm dy}{\sqrt{1-y^2}}+\sqrt{3}\left[\frac{y\arccos y}{\sqrt{1+3y^2}}\right]_{y=-\frac{1}3}^1 \\
    & = -\frac{1}2\left[\arcsin\Bigl(\frac{1-3y^2}2\Bigr)\right]_{y=-\frac{1}3}^1 +\frac{1}2\arccos\Bigl(-\frac{1}{3}\Bigr) =\frac{\pi}4 + \frac{1}2\arcsin\Bigl(\frac{1}3\Bigr)+\frac{1}2\arccos\Bigl(-\frac{1}{3}\Bigr).
\end{align}
Using the relation $\arccos + \arcsin = \pi/2$, we arrive at
$\mathcal W(\mu) = 12 \arccos\Bigl(-\frac{1}3\Bigr)$.
We now consider the points
\begin{align}
    x_{1}\vcentcolon= X\Big(\pi/3,-\arccos\Big(-\frac{\cot(\pi/3)}{\sqrt{3}}\Big)\Big) = \Big(\sqrt{2/3}, - 1/\sqrt{3},0\Big), \qquad x_2 = \big(-\sqrt{2/3}, -1/\sqrt{3},0\big).
\end{align}
Since the blow up at $x_1$ is given by the cone over the $1$-skeleton of a tetrahedron, it follows that the density is $\Theta^2(\mu,x_1) = \frac{3 \arccos{(-\frac{1}{3})}}{\pi}$, for instance using \Cref{lem:ap:Heppes} and \eqref{eq:theta_vs_length}. The same holds true for $x_2$ by symmetry.

\subsection{A branched immersion}\label{sec:branched_immersion}

\begin{figure}
     \centering
     \includegraphics[width=0.4\linewidth]{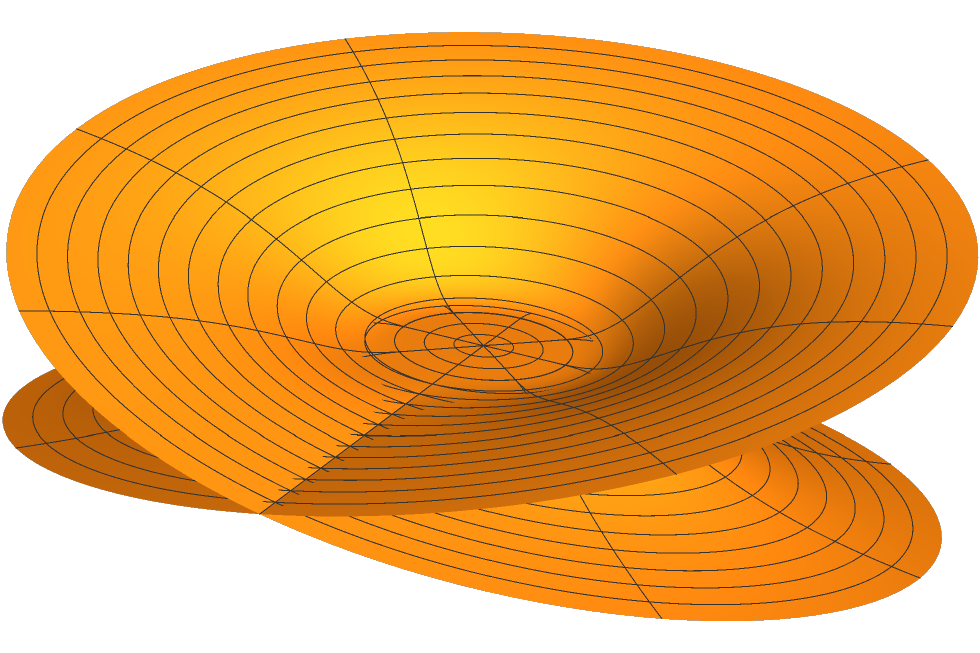}
     \caption{The branched immersion $f$ near the origin.}
     \label{fig:branched}
\end{figure}

Let $\delta \geq 0$, $\rho_0>0$, let $\psi\colon \R\to\R$ be a smooth cutoff function with $\psi(\rho)=1$ for $\rho\leq \rho_0/3$ and $\psi(\rho)\equiv 0$ for $\rho\geq 2\rho_0/3$. Then consider
\begin{align}
    f\colon [0,\infty)\times [0,2\pi)\to \RR^3, f(\rho,\theta) \vcentcolon = \big(\rho\cos 2\theta,\rho\sin 2\theta, \delta e^{-1/\rho^2}\psi(\rho)\cos{\theta}\big),
\end{align}
see \Cref{fig:branched}. 
By a short computation, $f$ is an immersion on $(0,\infty)\times [0,2\pi)$ since
\begin{align}
    |\partial_\rho f(\rho,\theta)\times \partial_\theta f(\rho,\theta)|^2\geq 4\rho^2>0.
\end{align}
Moreover, $f$ is Lipschitz and, by \cite[Example 5.9]{MR1704565}, the induced varifold $\mu$ is an integral curvature varifold with $B\in L^\infty(\mu)$ and $\Theta^2(\mu,0)=2$. 

Suppose that $F\in \mathcal{E}_\Sigma$ is some conformal Lipschitz immersion which parametrizes $\mu$ near the origin, i.e., $\spt\,\mu\cap B_\varepsilon(0) = F(\Sigma)\cap B_\varepsilon(0)$, $\varepsilon>0$. Let $U \vcentcolon= F^{-1}(B_\varepsilon(0))$ and let $p_1\neq p_2\in U$ with $F(p_1)=F(p_2)=0$. We have $H_F = H_{\mu}\circ F \in L^\infty(U)$ by \cite[Theorem 4.1]{MR2472179}. Now, the conformal transformation of the Laplace--Beltrami operator gives 
\begin{align}
e^{-2w} \Delta F= \Delta_{F^*\langle\cdot,\cdot\rangle}F = H_\mu \circ F,    
\end{align}
where $\Delta$ is the flat Laplacian in $\R^2$. Since $w\in L^\infty$, elliptic regularity implies that $F\in C^{1,\alpha}(U)$ for any $\alpha>0$, after shrinking $U$ if necessary. Hence, $F$ is an embedding in neighborhoods $U_1,U_2$ of $p_1,p_2$, respectively. Reducing $\varepsilon>0$ and shrinking $U_1,U_2$ we may assume that $\spt\,\mu \cap B_\varepsilon(0) = F(U_1)\cup F(U_2)$. We might therefore locally write $F(U_i)$ as a graph over its tangent space $T_0 F(U_i)$ for $i=1,2$. This implies that, for $\varepsilon>0$ small enough, $\spt\,\mu\cap \partial B_\varepsilon(0)$ is the union of two closed $C^1$-curves. However, since for $\delta>0$, $\spt\,\mu\cap \partial B_\varepsilon(0)$ is a single figure-eight-type curve, we obtain a contradiction.

The map $f$ gives rise to a global Lipschitz map $F\colon \mathbb{S}^2 \to \R^3$ as follows. We view $\mathbb{S}^2\subset\R^3$ with the north pole $N$ being the origin and $S=(0,0,-2)$. 
Flatten a tiny round disk of radius $\rho_0$ around $N$ (as in \cite[Section 6]{müller2023short}). Then, consider the map which doubles the polar angle (the angle in the plane $\{z=0\}$) which defines a smooth map and an immersion away from the $z$-axis. Hence, we obtain a smooth immersion on $\mathbb{S}^2\setminus \{N,S\}$, which can be checked to be Lipschitz on $\mathbb{S}^2$. Note that $f$ precisely realizes the composition of the flattening with this angle doubling in $B_{\rho_0}(0)\setminus B_{2\rho_0/3}(0)$. Thus, we may define a Lipschitz map $F\colon \mathbb{S}^2\to\R^3$ which is a smooth immersion on $\mathbb{S}^2\setminus\{N,S\}$ by replacing the flat disk with the image of $f$. Consider the Radon measure $\hat{\mu}\vcentcolon = F_\sharp \mathbb{S}^2$. Using the same arguments as in \cite[Example 5.9]{MR1704565} near $S$, we see that $\hat{\mu}$ is a integral curvature varifold in $\R^3$ with $B_{\hat\mu}\in L^\infty(\hat\mu)$. Moreover, its 
Willmore energy is arbitrarily close to $8\pi$ for $\delta,\rho_0>0$ small enough. However, by the discussion above, $\hat{\mu}$ may not be parameterized by any $W^{2,2}$-conformal Lipschitz immersion near the origin.



In order to obtain sharpness in \Cref{thm::main_8pi_curv}, one might also be tempted to consider the map $\mathbb C\cup\{\infty\}\to\mathbb C\cup\{\infty\}$ with $z\mapsto z^2$. Composed with the stereographic projection, this map induces a branched $C^\infty$-immersion on the Riemann sphere  $G\colon\mathbb S^2 \to \mathbb S^2\hookrightarrow \R^3$ with two branch points (one at each pole), $G(\mathbb S^2)=\mathbb S^2$, and $\mathcal W(G)=8\pi$. However, as opposed to the example described above, the varifold induced by $G$ can also be obtained by an unbranched $C^\infty$-embedding defined on $\Sigma$ given as the disjoint union of two copies of $\mathbb{S}^2$ or by equipping the standard embedding $\mathbb S^2\subset\R^3$ with multiplicity $2$. An honest example in higher codimension is given by the immersion $\mathbb C\to\mathbb C^2$ with $z\mapsto (z,z^2)$ with only one branch point, see \cite[Section 3.1]{MR3366862}.

\subsection{Singular set with arbitrary Hausdorff dimension}\label{sec:singular_set}

We define the regular set $\operatorname{reg}\mu$ as the set of points $x_0\in\spt\,\mu$ such that the conclusion of \cite{BiZhou22} (see \Cref{thm:BiZhou_Scaled}) holds in some ball $B_r(x_0)$, $r>0$, i.e., where $\spt\,\mu$ can be locally parametrized by a $W^{2,2}$-conformal Lipschitz immersion. The singular set is then $\operatorname{sing}\mu \vcentcolon= \spt\,\mu\setminus \operatorname{reg}\mu$.

The regularity theorem by Bi--Zhou \cite{BiZhou22} implies local regularity of integral $2$-varifolds with $H\in L^2$ near points $x_0$ with $\Theta^2(\mu,x_0)=1$. However, also if $\Theta^2(\mu,x)\equiv k$ for all $x\in B_r(x_0)$, we may apply the regularity result to the integral varifold $\frac{1}{k}\mu \mres B_r(x_0)$ to conclude $x_0\in \operatorname{reg}\mu$.
While the singular set is topologically small, see \Cref{rem:regularity}\eqref{item:sing_nowhere_dense} below, it can be large in Hausdorff dimension already for varifolds with energy close to the critical threshold of $8\pi$ in \Cref{thm::main_8pi_curv}. 

\begin{example}\label{ex:singular_set}
    Let $\varepsilon>0$ and $A\subset \{z\in\R^2\colon |z|<1\}$ be compact. Then there exists 
    an integral $2$-varifold $\mu$ in $\R^3$ satisfying  $\mathcal W(\mu)<8\pi +\varepsilon$, $\mu(\R^3)<\infty$, $\Theta^2(\mu,x)\in\{1,2\}$ for all $x\in\spt\,\mu$, and
    \begin{equation}
        \{x\in B_1(0)\colon \Theta^2(\mu,x)=2\} = A\times\{0\}.
    \end{equation}
\end{example}

\begin{remark}\label{rem:regularity}
    \leavevmode
    \begin{enumerate}
        \item\label{item:sing_nowhere_dense} Since the density $\Theta^2$ is upper-semi-continuous (cf.\ \Cref{lem:ap:monotonicity}), it is locally constant on an open dense subset of $\spt\,\mu$. In particular, the singular set of an integral $2$-varifold with finite Willmore energy is nowhere dense in $\spt\,\mu$. See \cite[Remark 8.1(1)]{Allard} for the analogous remark in Allard's regularity theory.
        \item  Let $\Sigma\vcentcolon=\spt\,\mu$ and consider 
    \begin{equation}
        S \vcentcolon=\partial \{x\in\Sigma\colon \Theta^2(\mu,x)=2\} \subset\operatorname{sing}\mu,
    \end{equation}
    where the boundary is taken in the topological space $\Sigma$. 
    Given any $0\le s<2$, there exists a compact set $A\subset\{z\in\R^2\colon |z|<1\}$ of Cantor type such that $\partial A = A$, and $\mathcal H^s(A)=\infty$, see \cite[2.10.29]{Federer}. In other words, since $A\times\{0\}\subset S$, $\operatorname{sing}\mu$ can be of any Hausdorff dimension less than 2. 
    \end{enumerate}
\end{remark}

\begin{proof}[Proof of \Cref{ex:singular_set}]
    Let $0\le\eta \in C^\infty(\R^2)$ such that $\eta(z)=1$ for $|z|\le1$ and $\eta(z)=0$ for $|z|\ge2$. By \cite[Theorem 3.14]{MS2.1} there exists a function $0\le u\in C^\infty(\R^2)$ such that 
    \begin{equation}
        A = \{z\in \R^2\colon u(z) = 0\}.
    \end{equation}
    Depending on $\delta > 0$, we define 
    \begin{equation}
        f_\delta\colon \{z\in\R^2\colon |z|\le 2\}\to \R^3,\qquad f_\delta(z)\vcentcolon=(z,\delta\eta(z)u(z)).
    \end{equation}
    By \cite[Lemma 6.1]{müller2023short}, we can glue the surface $f_\delta(B_2(0))$ into a large round sphere resulting in a smoothly immersed surface $\Sigma_\delta\subset\R^3$ with $\lim_{\delta\to0+}\mathcal W(\Sigma_\delta)=4\pi$. For $\delta$ small enough, we can achieve that $\Sigma_\delta$ is embedded and $\Sigma_\delta\cap B_1(0) = f_\delta(B_1(0))$. By $\hat{\Sigma}_\delta$, we denote the same construction with $u_\delta$ replaced by $-u_\delta$.
    Then, for $\delta$ small enough, the varifold 
    \begin{equation}
        \mu \vcentcolon= \theta \mathcal H^2\mres(\Sigma_\delta \cup \hat{\Sigma}_\delta)
    \end{equation}
    with
    \begin{equation}
        \theta(x) = 
        \begin{cases}
            1 &\text{ for $x \in (\Sigma_\delta \cup \hat{\Sigma}_\delta)\setminus(\Sigma_\delta \cap \hat{\Sigma}_\delta)$}\\
            2 &\text{ for $x \in (\Sigma_\delta \cap \hat{\Sigma}_\delta)$}
        \end{cases}
    \end{equation}
    has the desired property.
\end{proof}

\subsection{Minimizers with circular boundaries}\label{subsec:ex_circ_boundary}

We now present an example showing that the set of admissible boundary data in \Cref{thm:boundary_reg} is nonempty. To that end, we first compute the boundary integral in the case of a circle with an orthogonal prescribed conormal.

\begin{lemma}\label{lem:boundary_circle}
	Let $\gamma$ be a circle in a plane $P$ in $\RR^3$ and let $v$ be orthogonal to $P$. Then for all $x_0\in \RR^3$ we have
	\begin{align}
		\Big\vert\int_{\gamma}\frac{\langle x-x_0, v\rangle}{|x-x_0|^2}\mathrm{d}\mathcal{H}^1(x)\Big\vert \leq \pi.
	\end{align}
\end{lemma}

\begin{proof}
	Without loss of generality, $\gamma$ is the unit circle in the plane $\{z=0\}$ and $v\equiv e_3$. Let $x_0 = (a,b,c)\in \RR^3$ and assume without loss of generality $b=0$. Then we have
    $\langle x-x_0,v\rangle = - c$ for $x\in \gamma$,
	whereas for $x= (\cos t, \sin t, 0)\in \gamma$, $t\in [0,2\pi]$, we compute
	\begin{align}
		|x-x_0|^2 = (\cos t-a)^2 + \sin^2 t + c^2 = 1 - 2 a \cos t + a^2 + c^2.
	\end{align}
	Using the periodicity of the cosine we have
	\begin{align}
		\int_{\gamma} \frac{\langle x-x_0, v\rangle}{|x-x_0|^2}\mathrm{d}\mathcal{H}^1(x) &= -c \int_0^{2\pi} (1-2a \cos t +a^2+c^2)^{-1} \mathrm{d}t \\ &= - c \int_0^\pi \Big( (1-2a \cos t + a^2+c^2)^{-1} + (1+2a \cos t +a^2+c^2)^{-1}\Big)\mathrm{d}t \\
		& =\vcentcolon I(-a) + I(a).
	\end{align}
	We now apply the substitution $u=\tan t/2$ so that $\cos t = \frac{1-u^2}{1+u^2}$ and $\mathrm{d}t =  \frac{2 \mathrm{d}u}{1+u^2}$. Hence we find
	\begin{align}
		I(a) &= -c \int_0^\infty \Big(1 + 2a \frac{1-u^2}{1+u^2} + a^2+c^2\Big)^{-1} \frac{2}{1+u^2}\mathrm{d}u \\
		&= -2c \int_0^\infty \Big((1+a)^2 + c^2 + ((1-a)^2 + c^2)u^2\Big)^{-1}\mathrm{d}u \\
		&= -2c \frac{\arctan\Big[\sqrt{\frac{(1-a)^2+c^2}{(1+a)^2+c^2}} u\Big]}{\sqrt{(1+a)^2+c^2}\cdot\sqrt{(1-a)^2+c^2}} \Bigg\vert_{u=0}^{u=\infty} = \frac{-c \pi}{\sqrt{(1+a)^2+c^2}\cdot\sqrt{(1-a)^2+c^2}},
	\end{align}
	so that $I(a)=I(-a)$. The maximum modulus of this function is easily seen to be $\pi/2$ at $c^2=1-a^2$ and the claim follows.
\end{proof}

Now, consider, for example, two boundary curves $\gamma_1, \gamma_2$ that are parallel geodesics on a cylinder and sufficiently close. Then, \eqref{eq:LY_boundary_small} is satisfied and \Cref{thm:boundary_reg} yields the existence of smooth minimizers.

\appendix

\section{}

For the convenience of the reader, we have gathered essential background material in this section.

\begin{theorem}[{\cite[Theorem 1.1, rescaled version]{BiZhou22}}]\label{thm:BiZhou_Scaled}
     Let $\mu$ be an integral $2$-varifold in $\R^n$ with $\W(\mu)<\infty$. Let $x_0\in \Sigma=\spt\,\mu$. If
     \begin{align}\label{density condition}
     \mu(B_r(x_0))\leq (1+\varepsilon)\pi r^2, \qquad    \int_{B_r(x_0)}|H|^2\,\mathrm d\mu \leq \varepsilon^2,
     \end{align}
     for some $\varepsilon>0$ sufficiently small, 
     then there exists a bi-Lipschitz conformal parametrization $f\colon D\coloneqq B_1(0)\subset \R^2\to f(D)\subset \Sigma$ satisfying
     \begin{enumerate}
      \item  $B_{r(1-\psi)}(x_0)\cap\Sigma\subset f(D)$;
     \item For any $x,y\in D$,
     \begin{align}
     r(1-\psi)|x-y|\leq |f(x)-f(y)|\le r(1+\psi)|x-y|;
     \end{align}
     \item Let $g\coloneqq \mathrm df\otimes\mathrm df$, then there exist $w\in W^{1,2}(D)\cap L^\infty(D)$ such that
             \[g=e^{2w}(\mathrm dx\otimes \mathrm dx+\mathrm dy\otimes \mathrm dy)\]
            and $$\|w- \log r\|_{L^\infty(D)}+\|\nabla w\|_{L^2(D)} + \|\nabla^2 w\|_{L^1(D)}\le \psi;$$
     \item $f\in W^{2,2}(D,\mathbb{R}^n)$ and
     \[\|f-r \emph{i}\|_{W^{2,2}(D)}\leq r \psi ,\]
   where $\emph{i}\colon D\rightarrow \mathbb{R}^n$ is a standard isometric embedding.
     \end{enumerate}
     Here $\psi=\psi(\varepsilon)$ is a positive function such that $\lim_{\varepsilon\to 0}\psi(\varepsilon)=0$.
\end{theorem}

\begin{lemma}[\protect{Geodesic net on the $2$-sphere \cite{Heppes,Taylor,AllardAlmgren}}]\label{lem:ap:Heppes}
    Suppose $\gamma$ is an integral $1$-varifold in $\mathbb S^2$ and $\Theta^1(\gamma,p)<2$ for all $p\in\mathbb S^2$. If $\gamma$ is stationary with total length $\gamma(\mathbb S^2)<4\pi$, then either
    \begin{enumerate}
        \item $\gamma(\mathbb S^2)=2\pi$ and  $\gamma$ is a great circle; or
        \item $\gamma(\mathbb S^2)=3\pi$ and $\gamma$ is given by three half great circles having their endpoints in common; or
        \item\label{item:tetrahedron} $\gamma(\mathbb S^2)=6\arccos(-1/3)$ and $\gamma$ is given by six geodesic arcs forming a regular tetrahedron.
    \end{enumerate}
\end{lemma}

\begin{proof}
    Indeed, by Allard--Almgren \cite{AllardAlmgren} it is known that $\gamma$ consists of geodesic segments such that the sum of the tension forces acting on the junctions is zero. Such graphs on $\mathbb S^2\subset\R^3$ whose edges are geodesic segments and whose vertices satisfy the balancing condition have been classified by Heppes \cite{Heppes}. The full list is given in \cite[p.\,501]{Taylor}. In the following, we extend that list by adding the total length $\gamma(\mathbb S^2)$ and comparing it with the real number $4\pi$. Except for the elementary cases in (1)-(4), the edge lengths are taken from \cite{Taylor}.
    \begin{enumerate}
    \item Great circle, 
    \begin{equation}
        \gamma(\mathbb S^2) = 2\pi < 4\pi.
    \end{equation}
    \item Three half circles,
    \begin{equation}
        \gamma(\mathbb S^2) = 3\pi < 4\pi.
    \end{equation}
    \item Tetrahedron, six arcs,
    \begin{equation}
        \gamma(\mathbb S^2)=6\arccos\Bigl(-\frac{1}3\Bigr) < 11.5 < 4\pi.
    \end{equation}
    \item
    Cube, twelve arcs,
    \begin{equation}
        \gamma(\mathbb S^2)=12\arccos\Bigl(\frac{1}3\Bigr)>14.5>4\pi.
    \end{equation}
    \item
    Prism over a regular pentagon, 15 arcs,
    \begin{equation}
        \gamma(\mathbb S^2)=10\arccos\Bigl(\frac{\sqrt{5}}{3}\Bigr)+5\arccos\Bigl(\frac{3-5\sqrt{5}/3}{5-\sqrt{5}}\Bigr) > 16 > 4\pi.
    \end{equation}
    \item 
    Prism over a regular triangle, 9 arcs,
    \begin{equation}
        \gamma(\mathbb S^2)=6\arccos\Bigl(-\frac{1}3\Bigr)+3\arccos\Bigl(\frac{7}9\Bigr)>13.5>4\pi.
    \end{equation}
    \item
    Dodecahedron, 30 arcs,
    \begin{equation}
        \gamma(\mathbb S^2)=30\arccos\Bigl(1-\frac{8}{3(1+\sqrt{5})^2}\Bigr)>21>4\pi.
    \end{equation}
    \item 
    24 arcs forming 2 regular quadrilaterals and 8 equal pentagons, with each quadrilateral surrounded by 4 pentagons and each pentagon surrounded by 4 pentagons and one quadrilateral,
    \begin{align}
        \gamma(\mathbb S^2)&=8\cdot2\arcsin\Bigl(\frac{1}{\sqrt{3}}\Bigr) + 8\cdot2\arcsin\Bigl(\frac{\sqrt{2-\sqrt{2}}}{\sqrt{3}}\Bigr) \\
        &\qquad + 8\cdot2\arcsin\Bigl(\sqrt{\frac{(\sqrt[4]{2}-1)^2}{6} + \frac{(2-\sqrt{2})^2}{12}}\Bigr) >20>4\pi.
    \end{align}
    \item
    18 arcs forming 4 equal pentagons and 4 equal quadrilaterals, with each quadrilateral surrounded by 3 pentagons and 1 quadrilateral, and each pentagon by 3 quadrilaterals and 2 pentagons,
    \begin{equation}
        \gamma(\mathbb S^2)=\Bigl(6\cdot83.80167087^\circ + 8\cdot58.25684287^\circ+4\cdot13.55944752^\circ\Bigr)\frac{2\pi}{360^\circ}>17.5>4\pi.
    \end{equation}
    \item
    21 arcs forming 3 regular quadrilaterals and 6 equal pentagons, with each quadrilateral surrounded by 4 pentagons and each pentagon by 2 quadrilaterals and 3 pentagons,
    \begin{equation}
        \gamma(\mathbb S^2)=12\cdot 2\arcsin\Bigl(\frac{1}{\sqrt{3}}\Bigr) + 6\cdot 2\arcsin\Bigl(\sqrt{3-\frac{\sqrt{6}}{6}}\Bigr) + 3\cdot 2\arcsin\Bigl(\frac{\sqrt{3}-\sqrt{2}}{2\sqrt{3}}\Bigr)>25>4\pi.
    \end{equation}
    \end{enumerate}
\end{proof}

\begin{lemma}[{Monotonicity inequality \cite[Appendix A]{MR2119722}}]\label{lem:ap:monotonicity}
Suppose $\mu$ is an integral $2$-varifold in an open set $U\subset\R^n$ with locally square integrable mean curvature. Then 
\begin{align}\label{eq:Li-Yau}
    &\frac{\mu(B_r(x))}{\pi r^2}+\frac{1}{2\pi} \int_{B_r(x)}r^{-2}\langle H(y),y-x\rangle \,\mathrm d\mu(y)\\
    &\le \frac{\mu(B_s(x))}{\pi s^2} + \frac{1}{16\pi}\int_{B_s(x)}|H|^2\,\mathrm d\mu +\frac{1}{2\pi} \int_{B_s(x)}s^{-2}\langle H(y),y-x\rangle \,\mathrm d\mu(y) 
\end{align}
for all $x\in U$ and $0<r<s<\infty$ with $B_s(x)\subset U$. Moreover, the density $\Theta^2(\mu,x)$ exists at all $x\in U$ and satisfies
\begin{equation}
    \limsup_{y\to x}\Theta^2(\mu,y)\le \Theta^2(\mu,x),\qquad 1\le\Theta^2(\mu,x)<\infty \quad \text{for $x\in\spt\,\mu$}.
\end{equation}
If $U=\R^n$ and $\mu(\R^n)<\infty$, we may take $r\to0$ and $s\to\infty$ to obtain the \emph{Li--Yau-inequality}
\begin{equation}\label{eq:LY_closed}
    \Theta^2(\mu,x)\le\frac{1}{4\pi}\mathcal W(\mu).
\end{equation}
\end{lemma}

\section*{Acknowledgements}
This research was funded in whole, or in part, by the Austrian Science Fund (FWF), grant numbers \href{https://doi.org/10.55776/P32788}{10.55776/P32788} and \href{https://doi.org/10.55776/ESP557}{10.55776/ESP557}. 
The authors would like to thank Ernst Kuwert, Elena Mäder-Baumdicker, Ulrich Menne, and Matthias Röger for helpful comments and discussions.

\bibliography{Lib_new}
\bibliographystyle{abbrev}
\end{document}